\documentclass[12pt,a4paper]{article}
\usepackage[]{fontenc}
\usepackage[latin1]{inputenc}
\usepackage{a4,amsmath,amssymb,amsfonts,amsthm}
\usepackage{amsrefs}
\usepackage[all]{xy}
\usepackage{palatino}
 \usepackage[colorlinks]{hyperref}
 \usepackage{amsrefs}

\begin{document}
\bibliographystyle{plain}


\def\mR{\M{R}}
\def\mZ{\M{Z}}
\def\mN{\M{N}}           
\def\mQ{\M{Q}}
\def\mC{\M{C}}
\def\mG{\M{G}}



\def\Spec{{\rm Spec}}
\def\rg{{\rm rg}}
\def\Hom{{\rm Hom}}
\def\Aut{{\rm Aut}}
 \def\Tr{{\rm Tr}}
 \def\Exp{{\rm Exp}}
 \def\Gal{{\rm Gal}}
 \def\End{{\rm End}}
 \def\det{{{\rm det}}}
 \def\Td{{\rm Td}}
 \def\ch{{\rm ch}}
 \def\che{{\rm ch}_{\rm eq}}
  \def\Spec{{\rm Spec}}
\def\Id{{\rm Id}}
\def\Zar{{\rm Zar}}
\def\Supp{{\rm Supp}}
\def\eq{{\rm eq}}
\def\Ann{{\rm Ann}}
\def\LT{{\rm LT}}
\def\Pic{{\rm Pic}}
\def\rg{{\rm rg}}
\def\et{{\rm et}}
\def\sep{{\rm sep}}
\def\ppcm{{\rm ppcm}}
\def\ord{{\rm ord}}
\def\Gr{{\rm Gr}}
\def\ker{{\rm ker}}
\def\rk{{\rm rk}}
\def\ev{{\rm ev}}
\def\div{{\rm div}}


\def\beginProof{\par{\bf Proof. }}
 \def\endProof{${\qed}$\par\smallskip}
 \def\mtr#1{\overline{#1}}
 \def\ra{\rightarrow}
 \def\mfp{{\mathfrak p}}

 \def\mQ{{\Bbb Q}}
 \def\mR{{\Bbb R}}
 \def\mZ{{\Bbb Z}}
 \def\mC{{\Bbb C}}
 \def\mN{{\Bbb N}}
 \def\mF{{\Bbb F}}
 \def\mA{{\Bbb A}}
  \def\mG{{\Bbb G}}
 \def\CI{{\cal I}}
 \def\CE{{\cal E}}
 \def\CJ{{\cal J}}
 \def\CH{{\cal H}}
 \def\CO{{\cal O}}
 \def\CA{{\cal A}}
 \def\CB{{\cal B}}
 \def\CC{{\cal C}}
 \def\CK{{\cal K}}
 \def\CL{{\cal L}}
 \def\CI{{\cal I}}
 \def\CM{{\cal M}}
\def\CP{{\cal P}}
\def\CR{{\cal R}}
\def\CG{{\cal G}}
\def\CT{{\cal G}}
 \def\wt#1{\widetilde{#1}}
 \def\mod{{\rm mod\ }}
 \def\refeq#1{(\ref{#1})}
 \def\blb{{\big(}}
 \def\brb{{\big)}}
\def\mc{{{\mathfrak c}}}
\def\mcpr{{{\mathfrak c}'}}
\def\mcprpr{{{\mathfrak c}''}}
\def\ss{{\rm ss}}
\def\parf{{\rm parf}}
\def\P1{{{\bf P}^1}}
\def\cod{{\rm cod}}
\def\ss{\scriptstyle}
\def\OX{{ {\cal O}_X}}
\def\mpartial{{\mtr{\partial}}}
\def\inv{{\rm inv}}
\def\indlim{\underrightarrow{\lim}}
\def\Ind{{\rm Ind}}
\def\Ens{{\rm Ens}}
\def\without{\backslash}
\def\pbdb{{'o_b\ D^-_c}}
\def\qc{{\rm qc}}
\def\Com{{\rm Com}}
\def\an{{\rm an}}
\def\gfield{{\rm\bf k}}
\def\s{{\rm s}}
\def\dR{{\rm dR}}
\def\ari#1{\widehat{#1}}
\def\ul#1{\underline{#1}}
\def\sul#1{\underline{\scriptsize #1}}
\def\mou{{\mathfrak u}}
\def\ich{\mathfrak{ch}}
\def\cl{{\rm cl}}
\def\K{{\rm K}}
\def\R{{\rm R}}
\def\F{{\rm F}}
\def\L{{\rm L}}
\def\pgcd{{\rm pgcd}}
\def\rc{{\rm c}}
\def\N{{\rm N}}
\def\E{{\rm E}}
\def\H{{\rm H}}
\def\CHOW{{\rm CH}}
\def\A{{\rm A}}
\def\d{{\rm d}}
\def\Res{{\rm  Res}}
\def\GL{{\rm GL}}
\def\Alb{{\rm Alb}}
\def\alb{{\rm alb}}
\def\Hdg{{\rm Hdg}}
\def\Num{{\rm Num}}
\def\Irr{{\rm Irr}}
\def\Frac{{\rm Frac}}
\def\Sym{{\rm Sym}}
\def\indlim{\underrightarrow{\lim}}
\def\red{{\rm red}}
\def\naive{{\rm naive}}
\def\ch{{\rm ch}}
\def\Td{{\rm Td}}
\def\T{{\rm T}}
\def\min{{\rm min}}
\def\slope{{\rm slope}}
\def\char{{\rm char}}
\def\max{{\rm max}}
\def\min{{\rm min}}
\def\Sup{{\rm Sup}}
\def\Qb{\bar{\mQ}}
\def\mn{{\mu_n}}
\def\Zn{{\mZ/n}}
\def\bmn{\mn(\mC)}
\def\lmn{_\mn}
\def\dbd{{i\over 2\pi}\partial\bar\partial}
\def\umn{^\mn}
\def\op#1{\operatorname{#1}}
\def\prim{{\rm prim}}
\def\ach{{\ari{\rm ch}}}
\def\uZn{(\Zn)^\ast}
\def\Qmn{{\mQ_\mn}}
\def\ac1{\ari{\rm c}^1}
\def\odd{{\rm odd}}
\def\length{{\rm length}}
\def\Rat{{\rm Rat}}
\def\tot{{\rm tot}}
\def\Prim{{\rm Prim}}
\def\Mat{{\rm Mat}}


\def\RHom{{\rm RHom}}
\def\rRHom{{\mathcal RHom}}
\def\rHom{{\mathcal Hom}}
\def\dotimes{{\overline{\otimes}}}
\def\Ext{{\rm Ext}}
\def\rExt{{\mathcal Ext}}
\def\Tor{{\rm Tor}}
\def\rTor{{\mathcal Tor}}
\def\SP{{\mathfrak S}}
\def\perf{{\rm perf}}
\def\top{{\rm top}}
\def\Db{{\rm Dlb}}

\def\H{{\rm H}}
\def\D{{\rm D}}
\def\Del{{\mathfrak D}}
\def\Stab{{\rm Stab}}
\def\Div{{\rm Div}}
\def\Ver{{\rm Ver}}
\def\insep{{\rm insep}}
\def\Isog{{\rm Isog}}
\def\sgn{{\rm sign}}

 \newtheorem{theor}{Theorem}[section]
 \newtheorem{prop}[theor]{Proposition}
 \newtheorem{propdef}[theor]{Proposition-Definition}
 \newtheorem{cor}[theor]{Corollary}
 \newtheorem{lemma}[theor]{Lemma}
 \newtheorem{sublem}[theor]{sub-lemma}
 \newtheorem{defin}[theor]{Definition}
 \newtheorem{conj}[theor]{Conjecture}
 \newtheorem{rem}[theor]{Remark}
 \newtheorem{vconj}[theor]{Vague conjecture}
 \newtheorem{comp}[theor]{Complement}

 \parindent=0pt
 \parskip=5pt

 \author{Vincent MAILLOT \& Damian R\"OSSLER}
 \title{Conjectures on the logarithmic derivatives of Artin L-functions II}
\maketitle
\begin{abstract}
We formulate a general conjecture relating Chern classes of subbundles of Gauss-Manin bundles 
in Arakelov geometry to logarithmic derivatives of Artin $L$-functions of number fields. This conjecture may be viewed as a far-reaching generalisation of the \mbox{(Lerch-)Chowla}-Selberg formula computing logarithms of periods of elliptic curves in terms of special values of the $\Gamma$-function. We prove several special cases 
of this conjecture in the situation where the involved Artin characters are Dirichlet characters. 
This article contains the computations promised in \cite{Maillot-Roessler-Conjectures}, where 
our conjecture was announced. 
We also give a quick introduction to the Grothendieck-Riemann-Roch theorem and to the 
geometric fixed point formula, which form the geometric backbone of our conjecture.
\end{abstract}

\section{Introduction}

The main aim of this text is to provide the computations missing (and promised$\dots$) in the text \cite{Maillot-Roessler-Conjectures}. In that article, we formulated a conjecture, which relates the logarithmic derivatives of Artin $L$-functions at negative integers to certain Chern classes in Arakelov theory. This conjecture (see Conjectures \ref{RIconj} and 
Conjecture \ref{vRIconj} below) can be viewed as 
a far reaching generalisation of the (Lerch-)Chowla-Selberg formula, which 
computes the periods of CM elliptic curves in terms of special values of the $\Gamma$-function. 
A secondary aim of this article to provide 
a quick introduction to the main geometric ideas that lie behind our approach to Conjecture \ref{vRIconj}. These geometric ideas, although classical in some ways, are unfortunately not very well-known and we felt 
that to include a discussion of them here would make the computational part of the article (section 
\ref{seclogder}) more palatable. 

This text is an expanded version of some notes prepared by the second 
author for lectures given during two instructional conferences: the conference 
'Advanced Courses on Arakelov Geometry and Shimura Varieties' which took place at the Centre der Recerca Matem\`atica in Barcelona 
in February 2006 and the summer school 'Motives and complex multiplication', which took 
place in Ascona (Switzerland) in August 2016.

Very loosely speaking, the conjecture made in \cite{Maillot-Roessler-Conjectures} says the following. Suppose that you 
are provided with a relative pure homogenous polarised semistable (or log smooth) motive $\CM$ over an arithmetic base $B$ (which may have large dimension). Suppose also that $\CM$ carries the action of a number field $K$, with 
some compatibility with the polarisation. 
Then the Hodge realisation of this motive is a vector bundle $H$ on $B$, which is endowed with a (possibly mildly singular) hermitian metric coming from the polarisation. Furthermore, the vector 
bundle $H$ comes with an orthogonal decomposition 
$$H\simeq \bigoplus_{\sigma\in\Gal(K|\mQ)}H_\sigma.$$ 
Call $\bar H_\sigma$ the vector bundle $H_\sigma$ together with its hermitian metric. Arakelov theory 
associates with each $\bar H_\sigma$ its arithmetic Chern character $\ari{\ch}(\bar H_\sigma)$, which lives in the arithmetic Chow group $\ari{\CHOW}^\bullet_{\bar\mQ}(B)$ of $B$. 

Let now $\chi:\Gal(K|\mQ)\to\mC$ be an irreducible 
Artin character and let $l\geq 1$. 

{\it {\rm Conjecture}: the quantity} 
$$\sum_{\sigma \in\Gal(K|\mQ)}\ari{\ch}^{[l]}(\bar H_\sigma)\bar\chi(\sigma)$$
{\it is equal to the quantity } 
$$
2\,{L'(\chi,1-l)\over L(\chi,1-l)}+(1+{1\over 2}+\cdots+{1\over l-1})
$$
{\it multiplied by a certain $\bar\mQ$-linear combination of ordinary Chern classes of subbundles of 
$H$.}

Here $\ari{\ch}^{[l]}(\cdot)$ is the degree $l$ part of the arithmetic Chern character. See Conjecture \ref{vRIconj} below 
for a slightly more technical (but still vague) formulation. For abelian schemes, 
we can make a completely precise conjecture: this is Conjecture \ref{RIconj}. It should be 
possible to make a precise conjecture for semiabelian but generically abelian schemes but 
this seems difficult to do at the present time for lack of a sufficiently general 
theory of arithmetic Chern classes for singular hermitian metrics. This theory 
is being built in the articles \cite{BKK} and \cite{BKK2}. 

For abelian varieties with complex multiplication by a CM field, the quantities 
 $\ari{\ch}^{[1]}(\bar H_\sigma)$ can be computed in terms of periods. Thus in this case the 
equality above computes some linear combination of logarithms of periods in terms 
of the logarithmic derivatives ${L'(\chi,0)\over L(\chi,0)}$ of the irreducible Artin characters $\chi$ of 
the CM field. When the abelian variety is an elliptic curve, one recovers (a slight variant of) 
the formula of (Lerch-)Chowla-Selberg (see \cite{Selberg-Chowla-On-Epstein}). If 
$l=1$ and $K$ is an abelian extension of $\mQ$ one recovers a variant of the 
period conjecture of Gross-Deligne \cite[p. 205]{Gross-Periods} (not to be confused with 
the conjecture of Deligne \cite{Deligne-Valeurs} relatings periods and values of $L$-functions). 
For $l>1$, the invariant $\ari{\ch}^{[l]}(\bar H_\sigma)$ cannot be interpreted in terms of 
classical invariants anymore. Section \ref{secex} collects examples of 
computations in the literature, which fall in the framework of our conjecture (up to some finite 
factors which depend on the choices of models).

\begin{rem}\rm It is important to see that our conjecture falls outside the grid of the conjectures 
of Beilinson, Deligne, Stark, Gross and others (see eg \cite{RSS-Beilinson}) on the values of $L$-functions of motives. 
This can be seen from the fact that we are concerned here with 
the quotient between the second and the first coefficient of the Taylor series of 
an $L$-function at a non negative integer. This quotient in particular concerns 
the second coefficient of the Taylor series of the $L$-function at a non negative integer, about which the conjectures of 
Beilinson and Deligne do not make any prediction. The case of CM abelian varieties is somewhat confusing in this context, because in this case (as explained above), the conjecture computes some linear combinations of logarithms of periods. On the other hand, periods  
appear in Deligne's conjecture (see \cite{Deligne-Valeurs}) relating the values of the $L$-function of a motive to 
its periods. For CM abelian varieties, this conjecture was proven by Blasius (see \cite{Blasius-Critical} and the bibliography therein): the $L$-function of a 
CM abelian variety is a Hecke $L$-function of the CM field and the values of this Hecke $L$-function 
can be related to the periods of the abelian variety. Hecke $L$-functions are 
very different from Artin $L$-functions though, as is witnessed by the fact that in this case 
our conjecture relates the logarithmic derivatives of 
Artin $L$-functions to the {\it logarithms} of the periods of the abelian variety, whereas the result of Blasius relates values of Hecke $L$-functions to 
the the periods themselves (not their logarithms).\end{rem}

Our main contribution in this paper is a proof of a stabilised form of our conjecture in the 
situation where the number field $K$ is an abelian extension of $\mQ$ (so that the Artin characters become Dirichlet characters) and where the motive is smooth and arises 
from a finite group action on a Gauss-Manin bundle of geometric origin.
See Theorem \ref{AGBF} below. We also prove a stronger form of 
the conjecture in the situation where the number field $K$ is an abelian extension of $\mQ$ and 
the motive is the motive of an abelian scheme. See Theorem \ref{RIpart} below. In both cases, we derive our 
results from the equivariant Grothendieck-Riemann-Roch theorem in Arakelov geometry, applied to 
the relative de Rham complex. This theorem was first proven in degree one in \cite{lrr1} 
and in full generality in \cite{Tang-Concentration} and \cite{GRS-Ar} (put together). More details on the history of this theorem 
(whose main contributors are Bismut, Gillet, Soul\'e and Faltings) are given in section \ref{secERRTA}.

The structure of the text is as follows. Sections \ref{secGRR} to \ref{secERRTA} do not contain any original material and have been included for pedagogical reasons. In section \ref{secGRR}, we give a very quick introduction to the 
Grothendieck-Riemann-Roch formula. This theorem, although quite famous, is not as well-known 
as it should be and is rarely part of a standard course on algebraic geometry. In section \ref{secthofix}, we explain the content 
of Thomason's geometric fixed point formula for the action of a diagonalisable group. 
This formula (and its forerunners) is also a central result of algebraic geometry, which 
is not widely known. These two theorems can be formally combined to 
obtain an equivariant extension of the Grothendieck-Riemann-Roch formula, which we formulate 
in section \ref{secERRT}. We also examine there what statement one obtains when 
this theorem is applied to the relative de Rham complex. The resulting statement is 
a relative equivariant form of the Gauss-Bonnet formula. This statement is the geometric 
heart of our approach to Conjecture \ref{vRIconj}. 
In section \ref{secAra}, we give a historical snapshot of Arakelov theory and in section \ref{secERRTA}, 
we give a precise formulation of the equivariant arithmetic Grothendieck-Riemann-Roch formula, 
ie the equivariant Grothendieck-Riemann-Roch formula in Arakelov geometry. This formula 
will be our central tool. 
In section \ref{seclogder}, we apply this formula to the relative de Rham complex and we obtain a lifting to Arakelov theory of the relative equivariant form of the Gauss-Bonnet formula (see 
equation \refeq{AGB}); applying finite Fourier theory and some elementary results of analytic number theory, we transform this formula in an 
equality between a linear combination of logarithmic derivatives of Dirichlet $L$-functions evaluated at negative integers on the one hand and a linear combination of 
Chern classes of subbundles of Gauss-Manin bundles on the other hand. 
The final formula naturally suggests the general conjecture \ref{RIconj}, which 
is also included in section \ref{seclogder}. In section \ref{secex}, we examine various 
results on logarithmic derivatives of $L$-functions that have appeared in the literature and we show 
that they are all compatible with our general conjecture. We also explain there what part of these results are an actually consequence of \refeq{AGB}. 

{\bf Ackowledgments.} We are very grateful to J. Fr\'esan, one of the editors of this volume, for having the patience to 
wait for the completion of the present text. The second author also had many interesting conversations with him about the conjectures presented here. We would also like to thank J.-B. Bost, H. Gillet and C. Soulé for their support over the years and G. Freixas y Montplet for his continuing interest. We also benefitted from J.-M. Bismut's and X. Ma's remarks.

\section{The Grothendieck-Riemann-Roch formula}

In this section, 'scheme' will be short for 'noetherian and separated scheme'. 

\label{secGRR}

Let $C$ be a smooth projective curve over $\mC$. Let 
$D:=\sum_i n_i D_i$ be a divisor on $C$. The simplest instance 
of the Grothendieck-Riemann-Roch formula is probably 
the well-known equality
\begin{equation}
\chi({\cal O}(D)):=\dim_\mC H^0(C,{\cal O}(D))-\dim_\mC H^1(C,{\cal O}(D))=
\deg\ D+1-g
\label{RRcurves}
\end{equation}
where $$\deg\ D:=\sum_i n_i$$ is the degree of $D$ and $$g:=\dim_\mC H^0(C,\Omega_C)$$ is 
the genus of $C$. One can show that 
$$
\deg\ D=\int_{C(\mC)}c^1({\cal O}(D))
$$
where 
$c^1({\cal O}(D))$ is the first Chern class of $D$ (see eg \cite[Appendix A.3]{Har}). Thus \refeq{RRcurves} can be construed as a formula 
for the Euler characteristic $\chi({\cal O}(D))$ in terms of integrals of cohomology classes.

The 
Grothendieck-Riemann-Roch formula provides a similar formula for the Euler characteristic 
of any vector bundle, on any  scheme satisfying certain conditions and in 
a relative setting. Furthermore, the Grothendieck-Riemann-Roch formula is 
universal in the sense that it is independent of the cohomology theory. 
The remainder of this section is dedicated to the formulation of this theorem (in a slightly restricted 
setting).

First a definition. 

\begin{defin}
Let $X$ be a scheme. The group $K_0(X)$ (resp. ${K_0'}(X)$) is the free abelian 
group generated by the isomorphism classes of coherent locally free sheaves
(resp. coherent sheaves) on $X$, with relations $E=E'+E''$ if 
there is a short exact sequence
$$
0\ra E'\ra E\ra E''\ra 0.
$$
\label{defK}
\end{defin}

We shall also call a coherent locally free sheaf a {\it vector bundle.}
 
The group $K_0(X)$ (resp. $K_0'(X)$) is called the Grothendieck group of coherent locally free sheaves (resp. coherent sheaves) on $X$. 
If $f:X\ra Y$ is a proper morphism of  schemes, we define the map 
of abelian groups $\R f_*:{K_0'}(X)\ra {K_0'}(Y)$ by 
the formula 
$$
\R f_*(E):=\sum_{k\geq 0}(-1)^k \R^k f_*(E).
$$
This is well-defined, because the existence of the long exact sequence 
in cohomology implies that we have $\R f_*E=\R f_*E'+\R f_*E''$ in ${K_0'}(Y)$, $E,E'$ and $E''$ are 
as in Definition \ref{defK}. The group $K_0(X)$ is a commutative ring under 
the tensor product $\otimes$ and 
$K_0'(X)$ has a natural $K_0(X)$-module structure. The obvious map $K_0(X)\ra {K_0'}(X)$
is an isomorphism if $X$ is regular (see \cite[Th. I.9]{Manin-K-functor} if $X$ carries and ample line bundle and \cite[Lemme 3.3]{Tho} for the general case). Via this isomorphism, 
we obtain a map $\R f_*:K_0(X)\ra K_0(Y)$, if both $X$ and $Y$ are regular. 
For any morphism $f:X\ra Y$ of schemes, there is a pull-back map $Lf^*:K_0(Y)\ra K_0(X)$, defined 
in the obvious way, which is a map of rings. 

A theory kindred to $K_0$-theory is Chow theory. We first need a definition.

\begin{defin}
Let $R$ be a one-dimensional domain. Let $K:=\Frac(R)$ and let 
$f\in K^\ast$. Define the order of $f$ by the formula
$$
\ord_R(f):=\length_R(R/aR)-\length_R(R/bR)
$$
where $a\in R$, $b\in R^\ast$ are such that $f=a/b$. 
\end{defin}
One can show that the definition of $\ord_R(f)$ does not depend on the choice of $a,b$. 
Here the symbol $\length_R(\cdot)$ refers to the length of an $R$-module. 
See \cite[Appendix A.1\,\&\,A.3]{Fulton} for more details.

 Suppose for the time of the present paragraph that $X$ is an integral scheme. If $f\in\kappa(X)^*$ is a non zero rational function on $X$, we may define a formal $\mZ$-linear combination of codimension one closed integral subschemes of $X$ by the formula
$$
\div(f):=\sum_{x\in X,\,\cod(\bar x)=1}\ord_{\CO_x}(f)\bar x.
$$

For $p\geq 0$, we let 
$Z^p(X)$ be the free abelian group on all integral closed subschemes of codimension $p$
of $X$. An element of $Z^p(X)$ is called a $p$-cycle. We let $\Rat^p(X)\subseteq Z^p(X)$ be the subgroup of elements of the form $\div(f)$, where $f\in k(Z)^*$ is a rational function on 
a closed integral subscheme $Z$ of codimension $p-1$ of $X$.

\begin{defin}
$\CHOW^p(X):=Z^p(X)/\Rat^p(X).$ 
\end{defin}
The group $\CHOW^p(X)$ is called the $p$-th {\it Chow group} of $X$ and 
we shall write $$\CHOW^\bullet(X):=\bigoplus_{p\geq 0}\CHOW^p(X).$$

If $V$ is a closed subscheme of $X$, we write $[V]$ for the cycle 
$$\sum_{v\in V,\,\bar v\,\textrm{irred. comp. of $V$}}\length_{\CO_{V,v}}(\CO_{V,v})\bar v$$ 
in $X.$ Here $v$ runs through the generic points of the irreducible components of $V$ and 
$\bar v$ denotes the Zariski closure of $v$. 

By work of Gillet and Soulé, if $X$ is also regular, the group $\CHOW^\bullet(X)_\mQ:=
\CHOW^\bullet(X)\otimes\mQ$ can be made into a commutative $\mN$-graded ring. 
If we denote by $\cdot$ the multiplication in this ring, then we have 
$[W]\cdot [Z]=[Y\cap Z]$, if $W$ and $Z$ are closed integral 
subschemes of $X$ intersecting transversally (see \cite[I.2]{SABK} for more details and references). 

If $f:X\ra Y$ is a proper 
morphism of schemes, there is a unique push-forward map $f_*:\CHOW^\bullet(X)\ra\CHOW^\bullet(Y)$ such that 
$$f_*([Z])=[\kappa(Z):\kappa(f_*(Z))]\cdot [f_*Z]$$ 
if $Z$ is a closed integral subscheme $Z$ of $X$
such that $\dim(f_*Z)=\dim(Z)$ and such that 
$$f_*([Z])=0$$ otherwise. See \cite[Ex. 20.1.3, p. 396]{Fulton} for details.
If $f$ is a flat morphism, 
there is a pull-back map 
$$
f^*:\CHOW^\bullet(Y)\ra\CHOW^\bullet(X)
$$
such that $f^*([Z])=[f^*Z]$. Again see \cite[p. 394]{Fulton} for details.


Suppose now that $X$ regular. 
There is a unique ring homomorphism 
$$
\ch:K_0(X)\ra\CHOW^\bullet(X)_\mQ
$$
called the {\it Chern character}, with the following properties:

- $\ch(\cdot)$ is compatible with pull-back by flat morphisms;

- if $Z$ is an integral closed subscheme of codimension one of $X$, then 
$$\ch({\cal O}(Z))=\exp([Z]):=1+[Z]+{1\over 2!}[Z]\cdot[Z]+\cdots.$$

There is also a unique map
$$
\Td: K_0(X)\ra\CHOW^\bullet(X)_\mQ^*
$$
called the {\it Todd class}, with the following properties:

- $\Td(\cdot)$ is compatible with pull-back by flat morphisms;

- $\Td(x+x')=\Td(x)\cdot \Td(x')$;

- if $Z$ is an integral closed subscheme of codimension one of $X$, then 
$$
\Td({\cal O}(Z))={[Z]\over 1-\exp(-[Z])}.
$$

Finally there is a unique map $c:K_0(X)\ra\CHOW^\bullet(X)_\mQ^*$, called 
{\it the total Chern class}, such that

- $c(\cdot)$ is compatible with pull-back by flat morphisms;

- $c(x+x')=c(x)\cdot c(x')$;

- if $Z$ is an integral closed subscheme of codimension one of $X$, then $c({\cal O}(Z))=1+[Z]$.

The element $c^p(x):=c(x)^{[p]}(x)$ (where $(\cdot)^{[p]}$ takes the $p$-th graded part) is called the 
$p$-th Chern class of $x\in K_0(X)$. For a vector bundle $E/X$, we have 
$$
\ch(E)=1+c^1(E)+{1\over 2}(c^1(E)^2-2c^2(E))
+{1\over 6}(c^1(E)^3-3c^1(E)\cdot c^2(E)+3c^3(E))+\dots
$$
and 
$$
\Td(E)=1+{1\over 2}c^1(E)+{1\over 12}(c^1(E)^2+c^2(E))+{1\over 24}c^1(E)\cdot c^2(E)+\cdots
$$
We can now formulate the 
Grothendieck-Riemann-Roch theorem for smooth morphisms:

\begin{theor}
Let $X$, $Y$ be regular schemes. Let $f:X\ra Y$ be a smooth and strongly projective $S$-morphism. 
Then
$$
\ch(\R f_*(x))=f_*(\Td(\T f)\cdot\ch(x))
$$
for any $x\in K_0(X)$. 
\label{GRR}
\end{theor} 
Here the vector bundle $\T f:=\Omega_f^\vee$ is the dual of the sheaf of differentials of $f$. 
The vector bundle $\T f:=\Omega_f^\vee$ is also called the relative tangent bundle of $f$. 

{\bf Example}. Let $X:=C$ be a smooth and projective curve of genus $g$ over $\mC$,  as at the beginning of this section. 
Let $Y:=\Spec\ \mC$. Notice that $\CHOW^\bullet(Y)=\CHOW^0(Y)\simeq\mZ$ and that 
the Chern character of a vector bundle over $S$ is simply its rank under this identification. 
If we apply Theorem \ref{GRR} to $E:={\cal O}(D)$, we obtain, 
\begin{eqnarray*}
\ch(\R f_*({\cal O}(D)))&=&\chi({\cal O}(D))=
f_*((1+{1\over 2}c^1(\T C))(1+c^1({\cal O}(D))))\\
&=&
f_*(c^1({\cal O}(D))-{1\over 2}c^1(\Omega_C))=\deg(D)-{1\over 2}(2g-2)\\
&=&\deg(D)+1-g
\end{eqnarray*}
and thus we have recovered formula \refeq{RRcurves}.

The smoothness assumption on $f$ in Theorem \ref{GRR} can be relaxed. Suppose that 
$f$ is a strongly projective (but necessarily smooth) $S$-morphism. Then 
$f:X\ra Y$ has a factorisation $$f:X\stackrel{j}{\ra} {\bf P}^r\times Y\stackrel{\pi}{\ra} Y,$$ where 
$j$ is a closed immersion and $\pi$ is the natural projection. Theorem \ref{GRR}
still holds as stated if one replaces $\Td(\T f)$ by $j^*(\Td(\T\pi))\cdot\Td(N)^{-1}$, where 
$N$ is the normal bundle of the closed immersion $j$. 
The expression $(j^*\Td(\T\pi))\Td(N)^{-1}$ can be shown to be independent of the 
factorisation of $f$ into $j$ and $\pi$. See \cite[VIII, §2]{SGA6} for this. The fact that Theorem \ref{GRR} holds in this generality 
is a fundamental insight of Grothendieck; it shows that the theorem can be proved 
by reduction to the case of immersions and to the case of 
the structural morphism of ordinary projective space. 

The Riemann-Roch theorem for curves was discovered by B. Riemann and his student G. Roch in the middle of the nineteenth century. In the 1950s, F. Hirzebruch generalised the theorem to higher dimensional manifolds (but not to a relative situation). See his book \cite{Hirz} for this, where more historical references are given and the genesis of the Todd class is also described. The general relative case was first treated 
in the seminar \cite{SGA6} (see also \cite{Borel-Serre-RR}). The presentation of the Grothendieck-Riemann-Roch theorem 
given here follows W. Fulton's book \cite[chap. 15]{Fulton}. 

\section{Thomason's fixed point formula}

In this section, 'scheme' will be short for 'noetherian and separated scheme'. 
\label{secthofix}
We shall review a special case of Thomason's fixed point formula 
\cite[Th. 3.5]{Tho}. 

In the next paragraph, we give a list of definitions and basic results. 
These can found at the beginning of \cite{Tho}.

Let $X$ be a scheme.  
Let $\mn$ be the diagonalisable group 
scheme over $\Spec\ \mZ$ which corresponds to the finite group $\mZ/(n)$. 
Suppose that $X$ carries a $\mn$-action.  
We shall write $K^\mn_0(X)$ for the Grothendieck group 
of coherent locally free sheaves on $X$ which carry a $\mn$-equivariant 
structure. The definition of this group is completely similar to the definition of the 
ordinary Grothendieck group of locally free sheaves (see Definition \ref{defK}). 
Replacing locally free sheaves by coherent sheaves in the 
definition of $K^\mn_0(X)$ leads to the group $K^{'\mn}_0(X)$ and there is an obvious 
$K^\mn_0(X)$-module structure on the group $K^{'\mn}_0(X)$. 
If $X$ is regular, the natural morphism $K^\mn_0(X)\to K^{'\mn}_0(X)$ is an isomorphism 
(see \cite[Lemme 3.3]{Tho}). 
If the $\mn$-equivariant structure of $X$ is trivial, 
then the datum of a $\mn$-equivariant structure 
on a locally free sheaf $E$ on $X$ is equivalent 
to the datum of a $\mZ/(n)$-grading of $E$. 
For any $\mn$-equivariant locally free sheaf $E$ on $X$, 
we write $$\Lambda_{-1}(E):=\sum_{k=0}^{\rk(E)}(-1)^k\Lambda^k(E)\in K^\mn_0(X),$$
where $\Lambda^k(E)$ is the $k$-th exterior power of $E$, endowed with its natural 
$\mu_n$-equivariant structure. 
If $\Spec\ \mZ$ is endowed with the trivial $\mn$-structure, there is a unique isomorphism of rings 
${K_{0}^\mn(\Spec\ \mZ)}
\simeq{\mZ}[T]/(1-T^{n})$ with the following 
property: it maps the structure sheaf 
of $\Spec\, \mZ$ endowed with a homogenous $\mZ/(n)$-grading of weight one 
to $T$
and it maps 
any locally free sheaf carrying a trivial equivariant 
structure to the corresponding element of $K_{0}(\Spec\ \mZ).$
 
The functor of fixed points associated to $X$ is by definition 
the functor 
$${\rm\bf Schemes}/S\ra{\rm\bf Sets}$$ described by the rule 
$$
T\ \mapsto\ X(T)_{\mn(T)}.
$$
Here $X(T)_{\mn(T)}$ is the set of elements of $X(T)$ which 
are fixed under each element of $\mn(T)$. The functor 
of fixed points is representable by a scheme $X_\mn$ 
and the canonical morphism $X_\mn\ra X$ is a 
closed immersion (see \cite[VIII, 6.5 d]{SGA3-II}). 
Furthermore, if $X$ is regular then the scheme $X_\mn$ is regular (see \cite[Prop. 3.1]{Tho}).
We 
shall denote by $i$ the immersion $X_\mn\hookrightarrow X$. 
If $X$ is regular, we shall write $N^\vee$ for the dual of the normal sheaf  
of the closed immersion $i$. 
It is a locally free sheaf on $X_\mn$ and 
carries a natural $\mn$-equivariant structure. 
This structure corresponds to a $\mn$-grading, since 
$X_\mn$ carries the trivial $\mn$-equivariant structure and 
it can be shown that the weight $0$ term 
of this grading vanishes (see \cite[Prop. 3.1]{Tho}). 

Let $f:X\ra Y$ be a morphism between $\mn$-equivariant schemes which respects 
the $\mn$-actions. If $f$ is proper then the morphism $f$ induces 
a direct image map $$\R f_*:{K_0'}^\mn(X)\ra {K_0'}^\mn(Y),$$ which is 
a homomorphism of groups and is uniquely determined by the 
fact that  
$$
\R f_*(E):=\sum_{k\geq 0}(-1)^k \R^k f_*(E)
$$
 for 
any $\mn$-equivariant coherent sheaf $E$ on $X$. Here, as before,  
$\R^k f_*(E)$ refers to the $k$-th higher direct image sheaf
of $E$ under $f$; the sheaves $\R^k f_*(E)$ are coherent and carry a 
natural $\mn$-equivariant structure. If $X$ and $Y$ are regular, the direct image 
map $\R f_*$ induces a map $K_0^\mn(X)\ra K_0^\mn(Y)$ that we shall also 
denote by the symbol $\R f_*$.

The morphism $f$ also induces a pull-back map 
$$\L f^*:K_0^\mn(Y)\ra K_0^\mn(X);$$ this is a ring morphism which sends 
a $\mn$-equivariant locally free sheaf $E$ on $Y$ to 
the class of the locally free sheaf $f^*(E)$ on $X$, endowed with 
its natural $\mn$-equivariant structure. For any 
elements $z\in K_0^{'\mn}(X)$ and $w\in K_0^\mn(Y)$, the projection 
formula $$\R f_*(z\cdot Lf^*(w))=w\cdot \R f_*(z)$$ holds (provided $f$ is proper). This 
implies that the group homomorphism $\R f_*$ is 
a morphism of $K_0^\mn(\Spec\ \mZ)$-modules, if the group 
$K_0^\mn(X)$ (resp. $K_0^\mn(Y)$) is endowed with 
the $K_0^\mn(\Spec\ \mZ)$-module structure induced by the 
pull-back map $K_0^\mn(\Spec\ \mZ)\ra K_0^\mn(X)$ (resp. 
$K_0^\mn(\Spec\,\mZ)\ra K_0^\mn(Y)$).

Fix $\zeta_n\in\mC$ a primitive $n$-root of unity. In the following theorem, we shall view $\mQ(\mn)$ as a ${K_{0}^\mn(\Spec\ \mZ)}$-algebra 
via the homomorphism sending $T$ to $\zeta_n$.

\begin{theor}
Let $X,Y$ be schemes with $\mn$-actions and let $f:X\to Y$ be a morphism 
compatible with the $\mn$-actions. Suppose that $X$ and $Y$ are regular and that $f$ is proper. 
Suppose that the $\mn$-action on $Y$ is trivial. 
Then 

{\rm (1)} The element $\lambda_{-1}(N^\vee)$ is a unit 
in the ring \mbox{$K^\mn_0(X_\mn)\otimes_{K_0^\mn(\Spec\ \mZ)}{\mQ(\mn)}$.} 

{\rm (2)} For any 
 \mbox{element $x\in K_0^\mn(X)$,} the equality 
$$
\R f_*(x)=
\R f^{\mn}_*(\Lambda_{-1}(N^\vee)^{-1}\cdot
\L i^*(x))
$$
holds in $K^{\mn}_0(Y)\otimes_{K_0^\mn(\Spec\ \mZ)}{\mQ(\mn)}$.
\label{fixedp}
\end{theor}

Notice the formal analogy between Theorem \ref{fixedp} and Theorem \ref{GRR}: 
$Li^*$ takes the place of the Chern character and $\Lambda_{-1}(N^\vee)^{-1}$ takes 
the place of the Todd class.

{\bf Example}. Suppose that $Y=\Spec\,\mC$ and that $X_\mn$ is finite over $Y$. 
Let $g:X\to X$ be the automorphism corresponding to $\zeta_n\in\mn(\mC)$. 
Note that in this case, we have $K^{\mn}_0(Y)\simeq K^{\mn}_0(\Spec\ \mZ)$ via 
the natural pull-back map so that there is an isomorphism 
$$I:K^{\mn}_0(Y)\otimes_{K_0^\mn(\Spec\ \mZ)}{\mQ(\mn)}\simeq\mQ(\mn).$$
We leave it to the reader to verify that if $V$ is a $\mn$-equivariant vector space over 
$\mC$, then $I(V\otimes 1)=\op{Trace}(g_*:V\to V).$ If $x=E$, where 
$E$ is a $\mn$-equivariant vector bundle on $X$, Theorem \ref{fixedp} gives the equality
\begin{align*}
\sum_{k\geq 0}(-1)^k \op{Trace}&(g_*:H^k(Y,E)\to H^k(Y,E)) \\
&=
\sum_{y\in Y_\mn(\mC)}{\op{Trace}(E_y)\over 
\sum_{t=0}^{\rk(\Omega_{Y,y})}(-1)^t\op{Trace}(g_*:\Lambda^t(\Omega_{Y,y})\to \Lambda^t(\Omega_{Y,y}))}
\end{align*}
It is an exercise of linear algebra to show that 
$$
\sum_{t=0}^{\rk(\Omega_{Y,y})}(-1)^t\op{Trace}(g_*:\Lambda^t(\Omega_{Y,y})\to \Lambda^t(\Omega_{Y,y}))=
\det(\Id-g_*:\Omega_{Y,y}\to\Omega_{Y,y})
$$
so that
$$
\sum_{k\geq 0}(-1)^k\op{Trace}(g_*:H^k(Y,E)\to H^k(Y,E))=
\sum_{y\in Y_\mn(\mC)}{\op{Trace}(E_y)\over 
\det(\Id-g_*:\Omega_{Y,y}\to \Omega_{Y,y})}
$$
This formula is a special case of the so-called 'Woods Hole' fixed point formula (see \cite[Letter 2-3 August 1964]{GS-Corr}).

\section{An equivariant extension of the Grothendieck-Riemann-Roch theorem}

In this section, 'scheme' will be short for 'noetherian and separated scheme'. 

\label{secERRT}

If we formally combine the Grothendieck-Riemann-Roch theorem and 
the fixed point theorem of Thomason, we obtain the following 
theorem. Fix a primitive root of unity $\zeta_n\in\mC$. 

\begin{theor}
Let $X$ and $Y$ be regular schemes. Suppose that $X$ and $Y$  are 
equipped with a $\mn$-action. Suppose also that the $\mn$-structure of $Y$ is trivial.
Let $f:X\ra Y$ be a $\mn$-equivariant proper morphism and suppose that $f^\mn$ is smooth and strongly projective. Then for any $x\in K_0^\mn(X)$, the formula 
$$
\ch_\mn(\R f_*(x))=f_*(\ch_\mn(\Lambda_{-1}(N^\vee))^{-1}\Td(\T f^\mn)\ch_\mn(x))
$$
holds in $\CHOW^\bullet(Y)_{\mQ(\mn)}$. 
\label{EGRR}
\end{theor}

Here again, $N$ refers to the normal bundle of the immersion $X_\mn\ra X$. 
If $E$ is a $\mn$-equivariant vector bundle on $X$, writing $E_k$ for the 
$k$-th graded piece of the restriction of $E$ to $X_\mn$, we define 
$$
\ch_\mn(E):=\sum_{k\in\Zn}\zeta_n^k\cdot \ch(E_k)\in\CHOW^\bullet(X_\mn)_{\mQ(\mn)}.
$$
The element $\ch_\mn(E)$ is called the {\it equivariant Chern character} of $E$.

{\bf Example. The generalised Gauss-Bonnet formula.} Suppose that 
the assumptions of Theorem \ref{EGRR} hold and that in addition $f$ is smooth. We shall apply 
Theorem \ref{EGRR} to the image in $K_0^\mn(X)$ of the relative de Rham complex of $f$, ie to the element 
$$\Lambda_{-1}(\Omega_f)=1-\Omega_f+\Lambda^2(\Omega_f)-\Lambda^3(\Omega_f)+\dots+(-1)^{\rk(\Omega_f)}\Lambda^{\rk(\Omega_f)}(\Omega_f).$$  
Recall that we have exact sequence on $X_\mn$ 
$$
0\ra N^\vee\ra \Omega_X|_{X_\mn}\ra\Omega_{X_\mn}\ra 0.
$$
One can show that the symbol $\Lambda_{-1}(\cdot)$ is multiplicative 
 on short exact sequences of vector bundles (exercise! Use the splitting principle). In particular, we have
$$
\Lambda_{-1}(N^\vee)\cdot \Lambda_{-1}(\Omega_{X_\mn})=\Lambda_{-1}(\Omega_X|_{X_\mn})
$$
in $K_0^\mn(X_\mn)$. A last point is that for any vector bundle $E$, there is an identity 
of characteristic classes
$$
\ch(\Lambda_{-1}(E))\Td(E^\vee))=c^\top(E^\vee).
$$
This identity is called the Borel-Serre identity - see \cite[Example 3.2.4, 3.2.5]{Fulton} for a proof. With a view 
to simplifying the right hand side of the equality in Theorem \ref{EGRR}, we now compute
$$
\ch_\mn(\Lambda_{-1}(N^\vee))^{-1}
\Td(\T f^\mn)\ch_\mn(\Lambda_{-1}(\Omega_f|_{X_\mn}))=\ch(\Lambda_{-1}(\Omega_{f^\mn}))\Td(\Omega_{f^\mn}^\vee))=c^\top(\T f^{\mn}).
$$
Thus, by Theorem \ref{EAGRR} we have
\begin{equation}
\sum_{p,q\geq 0}(-1)^{p+q}\ch_\mn(\R f^p(\Lambda^q(\Omega_f)))=f^\mn_*(c^\top(\T f^{\mn}))
\label{GBfor}
\end{equation}
Formula \refeq{GBfor} is an equivariant extension of the Gauss-Bonnet formula (see eg 
\cite[chap. III, ex. 3.8]{Wells-Diff} for the non-equivariant formula in a cohomological setting) and 
we shall see further below that the lifting of formula \refeq{GBfor} to Arakelov theory carries deep arithmetic information.

Suppose that $Y=\Spec\,\mC$. 
Let $g:X\to X$ be the automorphism corresponding to a primitive $n$-th root of unity $\zeta_n\in\mn(\mC)$, as in the example given in the last section. Formula \refeq{GBfor} together with the existence of the Hodge decomposition gives the identities
\begin{eqnarray*}
&&\sum_{i,j\geq 0}(-1)^{i+j}\Tr(g_*:H^i(X,\Omega^j_X)\to H^i(X,\Omega^j_X))\\
&=&\sum_{k\geq 0}(-1)^k \Tr(g_*:H^k(X(\mC),\mC)\to H^k(X(\mC),\mC))=\int_{X_\mn}c^\top(\T f^{\mn})
\end{eqnarray*}
where $H^k(X(\mC),\mC)$ is the $k$-th singular cohomology group of $X(\mC)$ with coefficients in 
$\mC$. In particular, if $X_\mn$ consists of a finite set of points, we have 
\begin{equation}
\sum_{k\geq 0}(-1)^k \Tr(g_*:H^k(X(\mC),\mC)\to H^k(X(\mC),\mC))=\#X_\mn(\mC).
\label{toplef}
\end{equation}
Formula \refeq{toplef} is just the classical topological Lefschetz fixed point formula 
applied to $X(\mC)$ and the endomorphism $g$.

\section{Arakelov geometry}

\label{secAra}

Arakelov geometry is an extension of scheme-theoretic algebraic 
geometry, where one tries to treat the places at infinity of a number field (corresponding to 
the archimedean valuations) on the same footing as the finite ones. 
To be more precise, consider a scheme $S$ which is proper over $\Spec\ \mZ$ 
and generically smooth. 
For each prime $p\in\Spec\ \mZ$, we then obtain by base-change 
a scheme $S_{\mZ_p}$ on the spectrum of the ring of $p$-adic integers $\mZ_p$. 
The set $S(\mQ_p)$ is then endowed with the following natural notion of 
distance. Let $P,R\in S(\mQ_p)$; by the valuative criterion of properness, 
we can uniquely extend $P$ and $Q$ to elements $\widetilde{P},\widetilde{R}$ 
of $S(\mZ_p)$. We can then define a distance $d(P,R)$ by the formula 
$$
d(P,R):=p^{-\sup\{k\in\mZ|\widetilde{P}(\mod\ p^k)=\widetilde{Q}(\mod\ p^k)\}}.
$$
This distance arises naturally from the scheme structure of $S$. 
No such distance is available for the set $S(\mC)$ and the strategy of 
Arakelov geometry is to equip $S(\mC)$, as well as the 
vector bundles thereon with a hermitian metric in order to make up for that 
lack. The scheme $S$ together with a metric on $S(\mC)$ 
is then understood as a 'compactification' of $S$, in the sense 
that it is supposed to live on the 'compactification' of $\Spec\ \mZ$  
obtained by formally adding the archimedean valuation. 
The introduction of hermitian metrics, which are purely analytic 
data, implies that Arakelov will rely on a lot of analysis to define 
direct images, intersection numbers, Chern classes etc. 
Here is the beginning of a  list of extensions of classical scheme-theoretic objects 
that have been worked out in the literature:

\hskip-1cm\begin{tabular}{r|l}
$S$ & $S$ with a hermitian metric on $S(\mC)$\\
$E$ a vector bundle on $S$ & $E$ a vector bundle on $S$ with a 
  hermitian metric on $E(\mC)$\\
   cycle $Z$ on $S$ & a cycle $Z$ on $S$ with a Green current for $Z(\mC)$\\
the degree of a variety & the height of a variety over a number field\\
the determinant of cohomology & the determinant of cohomology equipped with its 
    Quillen metric\\
the Todd class of $f$ & the arithmetic Todd class of $f$ multiplied by ($1$-$R(\T f)$)\\
\vdots&\vdots
\end{tabular}

\smallskip
Here $f$ is the morphism $S\to\Spec\, \mZ$. 

Many theorems of classical algebraic geometry have been extended to Arakelov 
theory. In particular, there are analogs of the Hilbert-Samuel theorem 
(see \cite{GS8} and \cite{Abbes}), of the Nakai-Moishezon criterrion for ampleness (see 
\cite{Zhang}), of the Grothendieck-Riemann-Roch 
theorem (see \cite{GS8}) and finally there is an analog of the equivariant Grothendieck-Riemann-Roch theorem, 
whose description is of the main aims of this text.

 Arakelov geometry started officially in S. Arakelov's paper \cite{Ara}, who developped an intersection theory for surfaces in the compactified setting. 
G. Faltings (see \cite{Faltings}) then proved a Riemann-Roch 
theorem in the framework of Arakelov's theory. After that L. Szpiro and his students proved many other results 
in the Arakelov theory of surfaces. See \cite{Szpiro-DIH} and 
also Lang's book \cite{Lang-Ara} for this.
The theory was 
then vastly generalised by H. Gillet and C. Soulé, who defined compactified 
Chow rings, Grothendieck groups and characteristic classes in all dimensions (see \cite{GS2} and \cite{GS3}). 
For an introduction to Arakelov geometry, see the book \cite{SABK}.

\section{An equivariant Riemann-Roch theorem in Arakelov geometry}

\label{secERRTA}

The aim of this section is to formulate the analog in Arakelov 
geometry of Theorem \ref{EGRR}. With the exception of the 
relative equivariant analytic torsion form, we shall define precisely all the objects that we need  
but the presentation will be very compact and this section should not be used 
as a partial introduction to higher dimensional Arakelov theory. For this, we recommend 
reading the first few chapters of the book \cite{SABK}. 

Let $D$ be a regular arithmetic ring. By this we mean a regular, 
excellent, Noetherian domain, together 
with a finite set $\cal S$ of injective ring homomorphisms of $D\ra{\Bbb C}$, 
which is invariant under complex conjugation. We fix a primitive root of unity 
$\zeta_n\in\mC$. 

Ley $n\geq 1$. We shall call {\it equivariant arithmetic variety} 
an integral scheme $X$ of finite type 
over $\Spec\ D$,  endowed with a $\mn$-equivariant structure over $D$ and such that 
there is an ample $\mn$-equivariant line bundle on $X$. We also require 
the fibre of $X$ over the generic point of $D$ to be smooth.

We shall write  
 $X({\Bbb C})$ for the set of complex points of the variety 
$\coprod_{\sigma\in{\cal S}}X\times_{D,\sigma}{\Bbb C}$, 
which naturally carries the structure of a complex manifold. 
The groups $\bmn$ acts 
on $X({\Bbb C})$ by holomorphic automorphisms and we shall 
write $g$ for the automorphism corresponding to $\zeta_{n}$.
 As we have seen in section \ref{secthofix}, the fixed point scheme 
$X_{\mn}$ is regular and there are
 natural isomorphisms of complex manifolds 
 $X_{\mn}({\Bbb C})\simeq (X({\Bbb C}))_{g}$, where 
$(X({\Bbb C}))_{g}$ is the set of fixed points of $X(\mC)$ under the action of $\bmn$. 
Complex conjugation induces an antiholomorphic 
automorphism of $X({\Bbb C})$ and $X_{\mu_{n}}({\Bbb C})$, both of 
which we denote by $F_{\infty}$.  

If $M$ is a complex manifold, we shall write ${\frak
A}^{p,p}(M)$ for the set of smooth complex differential forms $\omega$ of type $(p,p)$ on 
a complex manifold $M$ and 
$$
\widetilde{\frak A}(M):=\bigoplus_{p\geq 0}({\frak A}^{p,p}(M)/({\rm Im}\,\partial +{\rm Im}\,\mtr{\partial})).
$$
The operator $\dbd$ induces a $\mC$-linear endomorphism of $\widetilde{\frak A}(M)$ and 
we shall write $H^{p,p}(M)$ for its kernel. The space $H^{p,p}(M)$ is part of the Aeppli cohomology 
of $M$ (see eg \cite[par. 2]{YY-BC} for the definition of Aeppli cohomology and its relation to other cohomology theories).

 We shall write ${\frak A}^{p,p}(X_{\mu_{n}})$ for the subspace of  ${\frak
A}^{p,p}(X_\mn(\mC))$ consisting of smooth complex differential forms $\omega$ of type $(p,p)$,
such that $F_{\infty}^{*}\omega=(-1)^{p}\omega$ and
$$
\widetilde{\frak A}(X_{\mu_{n}}):=\bigoplus_{p\geq 0}({\frak A}^{p,p}(X_\mn)/({\rm Im}\,\partial +{\rm Im}\,\mtr{\partial}))
$$
Similarly, we shall write $H^{p,p}(X_\mn)$ for the kernel of $\dbd$ in $\widetilde{\frak A}(X_{\mu_{n}})$. Note that $H^{p,p}(X_\mn)$ is a subspace of $H^{p,p}(X_\mn(\mC))$. 

A {\it hermitian equivariant sheaf (resp. vector bundle)} on $X$ is 
a coherent sheaf (resp. a vector bundle) $E$ on $X$, assumed 
locally free on $X({\Bbb C})$, 
equipped with a $\mu_{n}$-action which lifts the action of 
$\mu_{n}$ on $X$ and 
a hermitian metric $h$ on the vector bundle $E({\Bbb C})$, which is invariant under
$F_{\infty}$ and $\mn$. We shall write $(E,h)$ or 
 $\mtr{E}$ for 
an hermitian equivariant sheaf (resp. vector bundle). There is 
a natural ${\Bbb Z}/(n)$-grading 
$E|_{X_{\mn}}\simeq\oplus_{k\in{\Bbb Z}/(n)}E_{k}$ on  
the restriction of $E$ to $X_{\mn}$, whose terms are orthogonal, because 
of the assumed $g$-invariance of the metric. For $k\in{\Bbb Z}/(n)$, we write $\mtr{E}_{k}$ for $E_k$ endowed with the induced metric. 
We also often write $\mtr{E}_{\mn}$ for 
$\mtr{E}_{0}$.

If $\mtr{V}=(V,h_V)$ is a hermitian vector bundle on $X_\mn$ 
we write $\ch(\mtr{V})$ for the differential form $\Tr(\exp(\Omega_{h_V}))$. 
Here $\Omega_{h_V}$ is the curvature form associated with the unique connection on $V(\mC)$ 
whose matrix is given locally by $\partial H\cdot H^{-1}$, where $H$ is the matrix 
of functions representing $h_V$ locally. The differential form $\ch(\mtr{V})$ is both 
$\partial$- and $\bar\partial$-closed and its class in Bott-Chern cohomology 
represents the Chern character of $V(\mC)$ in the Bott-Chern cohomology of $X_\mn(\mC).$ 
Recall also that there is a natural map from 
Bott-Chern cohomology to Aeppli cohomology (see again \cite[par. 2]{YY-BC})
so that $\ch(\mtr{V})$ may also be viewed as a differential form 
representative for the Chern character of $V(\mC)$ in the Aeppli cohomology of $X_\mn(\mC).$ From the differential form 
$\ch(\mtr{V})$, using the fundamental theorem on symmetric functions, we may 
define differential form representatives  in Bott-Chern cohomology of other linear combinations of  Chern classes, like the Todd class 
$\Td(\mtr{V})$ or the total Chern class $c(\mtr{V})$. 

If $(E,h)$ is a hermitian equivariant sheaf, we write $\ch_g(\mtr{E})$ for the equivariant Chern character form 
$$\ch_g(\mtr{E}):=\ch_g((E({\Bbb C}),h)):=\sum_{k\in\mZ/(n)}\zeta_n^k\ch(\mtr{E}_k).$$  
The symbol $\Td_{g}(\mtr{E})$ refers to the differential form 
$$\Td(\mtr{E}\lmn)\Big(\sum_{i\geq 0}(-1)^{i}\ch_{g}(\Lambda^{i}(\mtr{E}))\Big)^{-1}.$$
If $${\cal E}:0\ra E'\ra E\ra E''\ra 0$$ is an exact sequence 
of equivariant sheaves (resp. vector bundles), we shall write $\mtr{\cal E}$ for the
sequence
$\cal E$ together with a datum of $\bmn$- and $F_{\infty}$- invariant hermitian metrics on
$E'({\Bbb C})$, 
$E({\Bbb C})$ and 
$E''({\Bbb C})$. With $\mtr{\cal E}$ and $\ch_g$ is associated an {\it equivariant 
Bott-Chern secondary class} $\widetilde\ch_{g}(\mtr{\cal E})\in\widetilde{\frak A}(X_{\mn})$, 
which satisfies
the equation $$\dbd\widetilde\ch_{g}(\mtr{\cal E})=\ch_g(\mtr{E}')+
\ch_g(\mtr{E}'')-\ch_g(\mtr{E}).$$ This class 
is functorial for any morphism of arithmetic varieties and vanishes 
if the sequence $\mtr{\cal E}$ splits isometrically. See \cite[par. 3.3]{lrr1} for all this. 

\begin{defin}
The arithmetic equivariant Grothendieck 
group $\ari{K}^{'\mn}_{0}(X)$ (resp. 
$\ari{K}^{\mu_{n}}_{0}(X)$) of $X$ is the 
free abelian group generated by the elements of $\widetilde{\frak A}(X_{\mu_{n}})$ and by 
the equivariant isometry classes of hermitian equivariant sheaves
 (resp. vector bundles), together with the 
relations
\begin{description}
\item[(a)] for every exact sequence $\mtr{\cal E}$ as above, we have 
$\widetilde\ch_{g}(\mtr{\cal E})=\mtr{E}'-\mtr{E}+\mtr{E}''$;
\item[(b)] if $\eta\in \widetilde{\frak A}(X_{\mu_{n}})$ is the sum in $\widetilde{\frak A}(X_{\mu_{n}})$ of 
two elements $\eta'$ and $\eta''$, then $\eta=\eta'+\eta''$ in 
$\ari K^{'\mn}_{0}(X)$ (resp. $\ari{K}^{\mu_{n}}_{0}(X)$).
\end{description}
\label{defAEGG}
\end{defin}
We shall now define a ring structure on $\ari{K}_{0}^{\mn}(X)$. Let 
$\mtr{V}, \mtr{V}'$ be hermitian equivariant vector bundles. Let $\eta,\eta'$ be 
elements of $\widetilde{\frak A}(X_{\mn})$. We define a \mbox{product $\cdot$}  on 
$\ari{K}_{0}^{\mn}(X)$ by the rules
$$
\mtr{V}\cdot\mtr{V}':=\mtr{V}\otimes\mtr{V}'  
$$
$$
\mtr{V}\cdot\eta=\eta\cdot\mtr{V}:=\ch_g(\mtr{V})\wedge\eta
$$
and 
$$
\eta\cdot\eta':=\dbd\eta\wedge\eta'
$$ and we extend it by linearity.  We omit the proof that 
it is well-defined (see \cite[par. 4]{lrr1} for this). 
Notice  that the definition of $\ari{K}_{0}^{'\mn}(X)$ (resp. $\ari{K}_{0}^{\mn}(X)$) implies that 
there is an exact sequence of abelian groups 
\begin{equation}
\widetilde{\frak A}(X_{\mn})\ra \ari{K}_{0}^{\mn'}(X)\ra K^{'\mn}_{0}(X)\ra 0
\label{ForgSeq}
\end{equation}
(resp. 
$$
\widetilde{\frak A}(X_{\mn})\ra \ari{K}_{0}\umn(X)\ra K\umn_{0}(X)\ra 0
\ \ \ ),
$$
where 
$K^{'\mn}_{0}(X)$ (resp. $K\umn_{0}(X)$) is the Grothendieck group of $\mn$-equivariant 
coherent sheaves (resp. locally free sheaves) considered in section \ref{secthofix}.
Notice finally that there is 
a map from $\ari{K}_{0}^{\mu_{n}'}(X)$ to the space of complex 
closed differential forms,
 which is defined by the formula 
$$\ch_g(\mtr{E}+\kappa):=\ch_g(\mtr{E})+\dbd\kappa$$
(where $\mtr{E}$ an hermitian equivariant sheaf and $\kappa\in
\widetilde{\frak A}(X_{\mn})$). This map is well-defined and 
we shall denote it by $\ch_{g}(\cdot)$ as well. 
We have as before: 
if $X$ is regular then the natural morphism  $\ari{K}^{\mu_{n}}_{0}(X)\ra
\ari K^{'\mn}_{0}(X)$ is an isomorphism. See \cite[Prop. 4.2]{lrr1} for this.

Now let $f:X\ra Y$ be an equivariant projective morphism of relative dimension $d$ over $D$ of equivariant regular arithmetic 
varieties. We suppose that $f$ is smooth over the generic point of $D$. We suppose that $X(\mC)$ is endowed with a Kähler fibration structure with respect to $f(\mC)$; this is a family of Kähler metrics on 
the fibers of $f(\mC):X(\mC)\ra Y(\mC)$, satisfying a supplementary condition that we do not have 
the room to detail here (see \cite[par. 1]{BK-Higher} for details). It is encoded in a real closed $(1,1)$-form $\omega_f$ on $X(\mC)$. In particular, the datum of $\omega_f$ induces a hermitian 
metric on the relative tangent bundle $\T f(\mC)$. 
We shall see an example of such a structure in the applications. We suppose that 
$\omega_f$ is $g$-invariant.  
Suppose also that the action of $\mn$ on $Y$ is trivial and finally suppose  
that there is a $\mn$-equivariant line bundle over $X$, which is
 very ample relatively to $f$.  
 
 Let now $\mtr{E}:=(E,h)$ be an equivariant hermitian sheaf on $X$ and suppose 
that $\R^k f_*(E)_\mC$ is locally free for all $k\geq 0$. Let $\eta\in\widetilde{\frak A}(X_{\mn})$.

We let $\R^\bullet f_{*}\mtr{E}:=\sum_{k\geq 0}(-1)^k R^k f_{*}\mtr{E}$ be the alternating sum of the higher direct image sheaves, endowed with their 
natural equivariant structures and $L_{2}$-metrics. 
For each $y\in Y(\mC)$, the $L^2$-metric on $\R^{i}f_{*}{E}(\mC)_y\simeq 
H^i_{\overline{\partial}}(X(\mC)_y,E(\mC)|_{X(\mC)_y})$ is 
defined by the formula
\begin{equation}
{1\over (2\pi)^{d}}\int_{X(\mC)_y}h(s,t){\omega_f^{d}}
\label{L2Def}
\end{equation}
where $s$ and $t$ are harmonic sections (i.e. in the kernel of the 
Kodaira Laplacian $\overline{\partial}\overline{\partial}^*+
\overline{\partial}^*\overline{\partial}$) of 
$\Lambda^i(T^{*(0,1)}X(\mC)_y)\otimes E(\mC)|_{X(\mC)_y}$. 
This definition is meaningful because by Hodge theory 
there is exactly one harmonic 
representative in each cohomology class.

Consider the rule, which associates the element 
$\R f_{*}^\bullet\mtr{E}-T_{g}(X,\mtr{E})$ of 
$\ari K^{'\mn}_{0}({Y})$ to $\mtr{E}$ and the element 
$\int_{X({\Bbb C})_{g}}\Td_{g}(\mtr{\T f}(\mC))\eta\in\ari K^{'\mn}_{0}({Y})$ to
$\eta$.  Here $T_{g}(\mtr{E})\in\widetilde{{\frak A}}(Y)$ is the 
equivariant analytic torsion form defined at the beginning of \cite{BM-Hol}. Its definition is too involved to 
be given in its entirety here but we shall define below its component of 
\mbox{degree $0$.} 

For the proof of the following proposition, see \cite[Th. 6.2]{Tang-Concentration}. 
\begin{prop}
The above rule extends to a well defined group homomorphism 
$f_{*}:\ari K^{'\mn}_{0}(X)\ra \ari K^{'\mn}_{0}(Y)$.
\label{WFProp}
\end{prop}

Now for the definition of the component of degree $0$ of $T_{g}(X,\mtr{E})$. 
Let $$\square_{q}^{E}:=\overline{\partial}\overline{\partial}^*+
\overline{\partial}^*\overline{\partial}$$ be the Kodaira Laplacian, which acts on the
  $C^\infty$-sections of the $C^\infty$-vector bundle\linebreak
\mbox{$\Lambda^{q}T^{*(0,1)}X(\mC)_y\otimes E(\mC)|_{X(\mC)_y}$} on $X(\mC)_y$. This space of
 sections is equipped with the $L^{2}$-metric as above and the
 operator $\square_q^{E(\mC)|_{X(\mC)_y}}$ is symmetric for that metric; we let
 \linebreak\mbox{${\rm Sp}(\square_q^{E(\mC)|_{X(\mC)_y}})\subseteq\mR$} be the set of eigenvalues
 of $\square_q^{E(\mC)|_{X(\mC)_y}}$ (which is discrete and bounded 
from below - see \cite[chap. 2, Prop. 2.36]{BGV}) and we let $\op{Eig}_{\,q}^{E(\mC)|_{X(\mC)_y}}(\lambda)$ be
 the eigenspace associated with an eigenvalue $\lambda$ 
(which is finite-dimensional - see \cite[chap. 2, Prop. 2.36]{BGV}). For $s\in\mC$ with $\Re(s)$ sufficiently large, we define
 $$
 Z(\mtr{E}|_{X(\mC)_y},g,s):=\sum_{q\geq 1}(-1)^{q+1}q\sum_{\lambda\in{\rm
 Sp}(\square_{q}^{E(\mC)|_{X(\mC)_y}})\backslash\{0\}}
 {\rm Tr}(g^*|_{\op{Eig}_{\,q}^{E(\mC)|_{X(\mC)_y}}(\lambda)})\lambda^{-s}.
 $$
 As a function of $s$, the function $Z(\mtr{E}|_{X(\mC)_y},g,s)$ has a meromorphic continuation to the whole complex plane,
 which is holomorphic around $0$. The degree $0$-part of the equivariant analytic torsion 
 form $T_{g}(\mtr{E})$ is then the complex number $Z'(\mtr{E}|_{X(\mC)_y},g,0)$. If $\R^k f_*(E)_\mC$ is locally free for all $k\geq 0$ 
 (which is our assumption)  then it can be shown that \linebreak $Z'(\mtr{E}|_{X(\mC)_y},g,0)$  
 is a $C^\infty$-function of $y$.  

We shall need the definition (due to Gillet and Soulé) of 'compactified' Chow theory. 
Let $X$ be a regular arithmetic variety over $D$. Let $p\geq 0$. 
We shall write $D^{p,p}(X)$ for the space of 
complex currents of type $p,p$ on $X(\mC)$ on which 
$F_\infty^*$ acts by multiplication by $(-1)^p$. 
Now let $A$ be  a subring of $\mC$ and suppose that $\mQ\subseteq A$. 
If $Z$ is a cycle of codimension $p$ with coefficients in $A$  on $X$ (in other words, a formal linear combination 
of integral closed subschemes of codimension $p$ with coefficients in $A$), a {\it Green current} $g_Z$ for $Z$ is 
an element of $D^{p,p}(X)$, which satisfies the equation
$$
{i\over 2\pi}\partial\mtr{\partial}g_Z+\delta_{Z(\mC)}=\omega_Z
$$
where $\omega_Z$ is a differential form and $\delta_{Z(\mC)}$ is the Dirac current associated with $Z(\mC)$. See the beginning of \cite{GS2} for this.
\begin{defin}
Let $p\geq 0$. The arithmetic Chow group $\ari{\CHOW}^p_A(X)$ is 
the $A$-vector space generated by the ordered pairs 
$(Z,g_Z)$, where $Z$ is a cycle of codimension $p$ with coefficients in $A$ on $X$ and $g_Z$ is a Green current for $Z$, with 
the relations
\begin{description}
\item[(a)] $\lambda\cdot(Z,g_Z)+(Z',g_{Z'})=(\lambda\cdot Z+Z',\lambda\cdot g_{Z}+g_{Z'})$; 
\item[(b)] $(\div(f),-\log|f|^2+\partial u+\mtr{\partial}v)=0$;
\end{description}
where $f$ is a non-zero rational function defined on 
a closed integral subscheme of codimension $p-1$ of $X$ and  
$u$ (resp. $v$) is a complex current of type $(p-2,p-1)$ (resp. 
$(p-1,p-2)$) such that $F_\infty^*(\partial u+\mtr{\partial}v)=(-1)^{p-1}(\partial u+\mtr{\partial}v)$. 
\label{arichowcoeff}
\end{defin}

We shall write $\ari{\CHOW}^\bullet_A(X):=\oplus_{p\geq 0}\ari{\CHOW}^p_A(X)$. 

\begin{rem}\rm The arithmetic Chow group with coefficients in $A$ defined in 
Definition \ref{arichowcoeff} is a formal variant of the arithmetic Chow group introduced 
by Gillet and Soul\'e in \cite{GS2}. In \cite{GS-Ar-an} they also consider a group, which is similar to 
ours in the case $A=\mR$ (but not identical). The properties of the arithmetic Chow group with coefficients in $A$
 listed below are similar to the properties of the arithmetic Chow group introduced in \cite{GS2} (with the same proofs going through verbatim) and we shall always refer to \cite{GS2} for properties of our group, although strictly speaking a different group is treated there.
 \end{rem}
 
The group 
$\ari{\CHOW}^\bullet_A(X)$ is equipped with a natural $A$-algebra structure, such that 
$$(Z,g_Z)\cdot(Z',g_{Z'})=(Z\cap Z',g_Z*g_{Z'})$$ 
if $Z,Z'$ are integral and  intersect transversally. Here the symbol $\ast$ refers to 
the star product, whose definition is too involved to be given here. See \cite[par. 2.1]{GS2} for this. A 
special case of the star product is described in the next example below. 
 If $f:X\ra Y$ is a projective and generically smooth morphism over $D$ between regular arithmetic varieties, there is a push-forward map $$f_*:\ari{\CHOW}^\bullet_A(X)\ra\ari{\CHOW}_A^\bullet(Y),$$ such that 
 $$f_*(Z,g_Z)=(\deg(Z/f_*(Z))f_*(Z),f(\mC)_*(g_Z))$$ for every integral closed subscheme $Z$ of $X$ and Green current 
 $g_Z$ of $Z$. Here we set $\deg(Z/f_*(Z))=[\kappa(Z):\kappa(f_*(Z))]$ if $\dim(f_*(Z))=\dim(Z)$ and 
 $\deg(Z/f_*(Z))=0$ otherwise. The expression $f(\mC)_*(g_Z)$ refers to the push-forward of currents. See \cite[par. 3.6]{GS2} for details. For any morphism $f:X\to Y$ over $D$ between regular arithmetic varieties, 
 there is a pull-back map $f^*:\ari{\CHOW}^\bullet_A(Y)\ra\ari{\CHOW}_A^\bullet(X)$, whose definition 
 presents the same difficulties as the definition of the ring structure on $\ari{\CHOW}_A^\bullet(\cdot)$. 
 See \cite[par. 4.4]{GS2} for details.   
 
 It is an easy exercise to show that the map of $A$-modules 
 $\mC\ra\ari{\CHOW}^1_A(\mZ)$, defined by the recipe $z\mapsto (0,z)$ is an isomorphism.
 
 If $X$ is a regular arithmetic variety, there is a unique ring morphism 
$$
\ari{\ch}:\ari{K}_0(X)\ra\ari{\CHOW}^\bullet_\mQ(X)
$$
called the arithmetic Chern character, such that 

- $\ari{\ch}$ is compatible with pull-backs by  $D$-morphisms;

- $\ari{\ch}(\eta)=(0,\eta)$ if $\eta\in \widetilde{\frak A}(X)$;

- if $\mtr{L}=(L,h)$ is a hermitian line bundle on $X$ and $s$ a rational section of $L$ then 
$$\ari{\ch}(\mtr{L})=\exp((\div\ s,-\log\ h(s,s))).$$
See the beginning of \cite{GS3} for this.

{\bf Example.} Suppose in this example that $X$ is regular and projective and flat of relative dimension $1$ over $D=\mZ$. 
 Suppose also that $Z$ and $Z'$ are two integral closed subschemes of codimension $1$ of $X$, 
 which  intersect transversally, are flat over 
 $\Spec\ \mZ$ and do not 
 intersect on the generic fiber. As $Z(\mC)$ (resp. $Z'(\mC)$) consists of one 
 point $P$ (resp. $P'$), the last condition just says that $P\not=P'$ in 
 $X(\mC)$.
 
 Now equip 
 ${\cal O}(Z)$ (resp. ${\cal O}(Z')$) with a conjugation invariant hermitian metric $h$ (resp. $h'$) and let 
 $s$ be a section of ${\cal O}(Z)$ (resp. $s'$ be a section of ${\cal O}(Z')$) vanishing 
 exactly on $Z$ (resp. $Z'$). In this case, we have 
\begin{align*}
 (Z,-\log\ h(s,s))\,\cdot&\,(Z',-\log\ h'(s',s'))\\
 &=(Z\cap Z',-\log\ h(s(P',P'))\delta_{Z(\mC)}-c^1(\mtr{{\cal O}(Z)})\log\ h'(s',s'))
\end{align*}
in  $\ari{\CHOW}^\bullet_\mQ(X)$ and hence, if $f$ is the morphism $X\ra\Spec\ \mZ$, 
\begin{align*}
 f_*(\ari{c}^1(\mtr{{\cal O}(Z)})\cdot\ari{c}^1(\mtr{{\cal O}(Z')}))&\\
 =
 \Big(0,(2\sum_{p\in f_*(Z\cap Z')}&\#\op{length}(Z_{\mF_p}\cap Z'_{\mF_p})\log\ p) \\
 &-\ \log\ h(s(P',P'))-\int_{X(\mC)}c^1(\mtr{{\cal O}(Z)})\log\ h'(s',s')\Big).
\end{align*}
From the arithmetic Chern character, using the fundamental theorem on symmetric functions, we may also define an arithmetic Todd class
\mbox{$\ari{\Td}:\ari{K}_0(X)\ra\ari{\CHOW}^\bullet_\mQ(X)^*$} and an arithmetic total Chern class 
\mbox{$\ari{c}:\ari{K}_0(X)\ra\ari{\CHOW}^\bullet_\mQ(X)^*$}.

If $\mtr{E}$ is an equivariant hermitian vector bundle 
on a regular equivariant arithmetic variety $X$, we define the {\it equivariant arithmetic Chern character} by the formula
$$
\ari{\ch}_\mn(\mtr{E})=\ari{\ch}_{\mn,\zeta}(\mtr{E}):=\sum_{k\in\Zn}\zeta_n^k\ari{\ch}(\mtr{E}_k)\in\ari{\CHOW}^\bullet_{\mQ(\mn)}(X_\mn)
$$
We write as before 
$\Lambda_{-1}(\mtr{E}):=\sum_{k=0}^{{\rm rk}(E)}(-1)^{k}\Lambda^{k}(\mtr{E})
\in\ari{K}^{\mu_{n}}_{0}(X)$, where $\Lambda^{k}(\mtr{E})$ is the 
$k$-th exterior power of $\mtr{E}$, endowed with its natural 
hermitian and equivariant structure.  

Finally, to formulate the equivariant Grothendieck-Riemann-Roch theorem in Arakelov geometry, 
we shall need the following exotic characteristic class. Let $X$ be a regular arithmetic variety. 

Recall that for any $z\in\mC$ with $|z|=1$, the Lerch zeta function $\zeta_L(z,s)$ is defined by the formula
$$
\zeta_L(z,s):=\sum_{k\geq 1}{z^k\over k^s}
$$
which is naturally defined for $\Re(s)>1$ and can be meromorphically continued to the whole plane.  
For $n$ any positive integer, define the $n$-th {\it harmonic number} $\CH_n$ by the formula
$$
\CH_0:=0
$$
and 
$$
\CH_n:=(1+{1\over 2}+\dots+{1\over n})
$$
when $n > 0$. For any $z\in\mC$ we now define the formal complex power series
$$
\widetilde{R}(z,x):=\sum_{k\geq 0}\Big(2\zeta_L'(z,-k)+\CH_k\cdot \zeta_L(z,-k)\Big){x^k\over k!}.
$$
(for those $z\in\mC$ where it makes sense) 
and 
$$
R(z,x):={1\over 2}(\widetilde{R}(z,x)-\widetilde{R}(\mtr{z},-x)).
$$
For any fixed $z\in\mC$, we identify $R(z,x)$ (resp. $\widetilde{R}(z,x)$) with the unique additive cohomology class it defines 
in Aeppli cohomology. 
For a $\mn(\mC)$-equivariant vector bundle $E$ on $X(\mC)$, where $X(\mC)$ is endowed with the trivial $\mn(\mC)$-equivariant structure, we now define the cohomology class $R_g(E)$ (resp. $\widetilde{R}_g(E)$) on 
$X(\mC)_g$ by the formula 
$$
R_g(E):=\sum_{u\in\Zn}R(\zeta_n^u,E_u).
$$
(resp. 
$$
\widetilde{R}_g(E):=\sum_{u\in\Zn}\widetilde{R}(\zeta_n^u,E_u)\,).
$$
The class $R_g(E)$ is often called the $R_g$-genus of $E$. 
Note that by construction we have
$$
\widetilde{R}_g(E)=\sum_{k\geq 0}\sum_{u\in\Zn}\Big(2\zeta_L'(\zeta_n^u,-k)+\CH_k\cdot\zeta_L(\zeta_n^u,-k)\Big){\rm ch}^{[k]}(E_u)
$$
(resp. 
\begin{align*}
R_g(E)=\sum_{k\geq 0}\sum_{u\in\Zn}\Big((\zeta_L'(\zeta_n^u,-k)&-(-1)^k\zeta_L'(\bar\zeta_n^u,-k)) \\
&+{1\over 2}\CH_k\cdot(\zeta_L(\zeta_n^u,-k)-(-1)^k
\zeta_L(\bar\zeta_n^u,-k))\Big){\rm ch}^{[k]}(E_u)\,).
\end{align*}
Let now again $f:X\to Y$ be an equivariant projective morphism over $D$ between regular equivariant arithmetic varieties. 
Suppose that there is an equivariant relatively ample line bundle on $X$ and that the equivariant structure of $Y$ is trivial. Suppose also that $f^\mn:X_\mn\to Y$ is smooth. 

Let ${N}={N}_{X/X_{\mu_{n}}}$ be the normal 
bundle of $X_{\mu_{n}}$ in $X$, which has a natural $\mn$-equivariant structure. 
The bundle $N(\mC)$ is by construction a quotient 
of the restriction to $X(\mC)_g$ of the relative tangent bundle $\T f(\mC)$ and we thus endow 
it with the corresponding quotient metric structure (which is $F_{\infty}$-invariant). We refer to 
the resulting $\mn$-equivariant hermitian vector bundle as $\mtr{N}=\mtr{N}_{X/X_\mn}$.
\begin{theor}[equivariant arithmetic Riemann-Roch theorem] Suppose that $f^\mn$ is smooth. Then 
the equality 
$$
\ari{\ch}_\mn(f_*(x))=f^{\mn}_*(\ari{\ch}_\mn(\Lambda_{-1}(\mtr{N}^\vee))^{-1}\Td(\mtr{\T f}^\mn)(1-R_g(\T f))\ari{\ch}_\mn(x))
$$
holds in $\ari{\CHOW}^\bullet_{\mQ(\mn)}(Y)$, for any $x\in\ari{K}_0^\mn(X)$. 
\label{EAGRR}
\end{theor}

\begin{rem}\rm If $f$ is smooth then $f^\mn$ is smooth. We leave the proof of this statement as an exercise for the reader.\end{rem}
Theorem \ref{EAGRR} results from a formal combination of the main results of 
\cite{Tang-Concentration} and \cite{GRS-Ar}. It is important to underline that the most difficult 
part of the proof is analytic in nature and is contained in J.-M. Bismut's article \cite{BM-Hol}. 
A proof of the degree one part of Theorem \ref{EAGRR} is given in \cite{lrr1}.

\section{Logarithmic derivatives of Dirichlet $L$-functions and arithmetic Chern classes of Gauss-Manin bundles}

\label{seclogder}

In this section, we shall apply Theorem \ref{EAGRR} to the relative de Rham complex 
of a smooth and projective morphism of regular equivariant arithmetic varieties (satisfying certain conditions) and interpret the resulting equality 
in terms of logarithmic derivatives of Dirichlet $L$-functions. 

As usual, fix a primitive $n$-th root of unity $\zeta_n\in\mC$. For convenience, we shall write $\zeta=\zeta_n$ in this section. If $\sigma\in\uZn$, we shall often write 
$
\sigma(\zeta)$ for $\zeta^\sigma.
$ 
If $\chi$ is a primitive Dirichlet character modulo $n$ (see eg 
\cite[chap. 4]{Washington-Cyc} for an introduction to Dirichlet characters) we shall write
$$
\tau(\chi)=\tau_\zeta(\chi):=\sum_{\sigma\in\uZn}\sigma(\zeta)\chi(\sigma)
$$
for the corresponding Gauss sum. 
We shall also write 
$$
L(\chi,s)=\sum_{n=0}^\infty{\chi(n)\over n^s}
$$
for the $L$-function associated with $\chi$. This function is defined for $\Re(s)>1$ but can be 
meromorphically continued to the whole complex plane. The resulting function is holomorphic 
everywhere if $\chi$ is not the trivial character. 

The following combinatorial lemmata will be needed in the proof.

\begin{lemma}
 \label{lemCHL}
 Let $M$ be complex projective manifold and
 let $E$ be a vector bundle on $M$ together with an  
automorphism
 $g:E\ra E$ of finite order (acting fiberwise). Let
 $\kappa$ be the class 
 $$
 \kappa:=\Td(E_0){\sum_{p\geq 0}(-1)^p p\cdot \ch_g(\Lambda^p(E^{\vee}))\over
 \sum_{p\geq 0}(-1)^p\ch_g(\Lambda^p(E_{\not=0}^{\vee}))}.
 $$
in the Aeppli cohomology of $M$. Then the equality
\[
\kappa^{[l + \op{rk}(E_{0})]} = - c^{\op{top}}(E_{0})
\sum_{z\in\mC}
\zeta_{\op{L}}(z,-l)\op{ch}^{[l]}((E^{\vee})_z)
\]
 holds for all $l\geq 0.$
\end{lemma}
Her $E_z$ is the largest subbundle of $E$ on which $g$ acts by multiplication by $z$.
\beginProof
See \cite[Lemma 3.1]{MR-Periods}. \endProof
Let $M$ be a complex projective manifold and
 let $(L,h_L)$ be an ample line bundle on $M$, endowed with a positive metric $h_L$. It is interesting (and it will be necessary later) to have an explicit formula for the $L_2$-metric carried by the vector spaces $H^p(M,\Omega^q_M)$ ($p,q\geq 0$), where 
the $L_2$-metric is computed using the Kähler metric coming from 
$c^1((L,h_L))$ and the metric on $\Omega^q_M$ is induced by $c^1((L,h_L))$. 
 
 Let us
 denote by $\omega\in H^{2}(M,\mC)$ the first Chern class of $L$ and
 for $k\leqslant \dim(M)$, let us write
 $P^{k}(M,\mC)\subseteq H^{k}(M,\mC)$ for the
 primitive cohomology associated to $\omega$; this is a Hodge substructure 
 of $H^{k}(M,\mC)$. Recall that for any $k\geq 0$,
 the primitive decomposition
 theorem of Lefschetz establishes an isomorphism
 $$
 H^k(M,\mC)\simeq\bigoplus_{r\geq
 \max(k-d,0)}\,\omega^{r}\wedge P^{k-2r}(M,\mC).
 $$
 Define the cohomological star operator $$*:H^{k}(M,\mC)\ra
 H^{2d-k}(M,\mC)$$ by the rule $$*\,\omega^{r}\wedge \phi :=
 i^{p-q}(-1)^{(p+q)(p+q+1)/2}{r!\over
 (d-p-q-r)!}\omega^{d-p-q-r}\wedge \phi$$ if $\phi$ is a primitive
 element of pure Hodge type $(p,q)$ and extend it by additivity. We can now
 define a pairing on $H^{k}(M,\mC)$ by the
 formula
 $$
 (\nu,\eta)_L:={1\over(2\pi)^d}\int_{M}\nu\wedge *\,\mtr{\eta}
 $$
 for any $\nu,\eta\in H^{k}(M,\mC)$. This pairing turns out to be a hermitian metric, which is sometimes called the 
 {\it Hodge metric}.  See \cite{Gr-Var} for all this. 
 \begin{lemma}
\label{HDRLem}
 The Hodge-de Rham isomorphism $$H^k(M,\mC)\simeq\bigoplus_{p+q=k}
 H^{q}(M,{\Omega}_M^p)$$ is an isometry if the right-hand side is endowed with 
 the Hodge metric and the left-hand side with the $L_2$-metric.
 \end{lemma}
\beginProof
See \cite[Lemma 2.7]{MR-Periods} .
\endProof
\begin{cor}
Let $h:M\to N$ be a projective and smooth morphism between quasi-projective complex manifolds. Let 
$g$ be a finite automorphism of $M$ over $N$ (ie $g$ acts fiberwise). Let 
\mbox{$h_0:M_g\to N$} be the induced morphism (which is smooth). The equality of 
characteristic classes
$$
\sum_{p,q}p\cdot (-1)^{p+q}\cdot \ch_{g}(\R^q h_*(\Omega^p_h))=
-\int_{M_g/N}c^\top(\T h_0)\Big[\sum_{l\geq 0}\sum_{z\in\mC}\zeta_L(z,-l)\ch^{[l]}((\Omega_h|_{M_g})_{z})\Big]
$$
in Aeppli cohomology holds. 
\label{corGlob}
\end{cor}
\beginProof
This is an immediate consequence of Lemma \ref{lemCHL}, the main result of 
\cite[Th. 2.12]{AS-Index2} and 
the fact that there is a natural map from Hodge cohomology (also called 
$\bar\partial$-cohomology) to Aeppli cohomology 
(see again \cite{YY-BC}, especially the diagram (2.1)). 
\endProof
\begin{lemma}
For any primitive Dirichlet character modulo $n$ and any $u\in\mZ$, 
the equality
$$
\sum_{\sigma\in\uZn}\sigma(\zeta^u)\chi(\sigma)=\bar\chi(u)\tau(\chi)
$$
holds.
\label{lemGS}
\end{lemma}
\beginProof
See \cite[chap. 4, lemma 4.7]{Washington-Cyc}.
\endProof

Lemma \ref{lemGS} implies the following two lemmata:

\begin{lemma}
Let $\mtr{V}$ be a $\mn$-equivariant hermitian vector bundle on a regular arithmetic 
variety $Z$. Suppose that $\mn$ acts trivially on $Z$. Then for any primitive character modulo $n$ we have 
\begin{eqnarray*}
\sum_{\sigma\in\uZn}\ach_{\mn,\sigma(\zeta)}(\mtr{V})\chi(\sigma)=\tau(\chi)\sum_{u\in\Zn}\ach(\mtr{V}_u)\bar\chi(u).
\end{eqnarray*}
\label{lemCO}
\end{lemma}
\begin{lemma}
For any primitive character modulo $n$ we have
$$
\sum_{\sigma\in\uZn}\zeta_L(\sigma(\zeta^u),s)\chi(\sigma)=\tau(\chi)\bar\chi(u)\L(\bar\chi,s).
$$
\label{lemLChi}
\end{lemma}

\begin{rem}\rm({\bf Important}). 
Let us call the number $\chi(-1)\in\{1,-1\}$ the parity of the \mbox{Dirichlet} character $\chi$ and for
any (positive) integer $l\in\mZ$ let us define the parity of $l$ to be $(-1)^l\in\{1,-1\}$.
By classical results of analytic number theory, 
we have $L(\chi_{\rm Prim},1-l)\not=0$ if $\chi$ and $l$ have the same parity (see \cite[before Th. 4.2]{Washington-Cyc}). More generally if $\chi$ is now an Artin character attached to any finite dimensional complex irreducible representation $R_{\chi}$ of the Galois group of a finite Galois extension of $\mQ$; we will say that $\chi$ is even (resp. odd) if $R_{\chi}(F_\infty) =\Id$ (resp. $R_{\chi}(F_\infty) = -\Id$), where $F_{\infty}$ is acting as the complex conjugation. 
Let's then denote by $L(\chi,s)$ the Artin $L$-function associated with $\chi$ 
(cf. \cite{Tate-Stark} or \cite[\S\ 7.10-12]{Neu} for an introduction). The function $L(\chi,s)$ is nonvanishing for $\Re(s) > 1$ and by Brauer admits a functional equation and a meromorphic continuation to the whole complex plane. One easily deduces from this the zeroes of $L(\chi,s)$ lying on the real negative line 
(cf. for instance \cite[p.541]{Neu}). We get again that $L(\chi,1-l)\not=0$ when $\chi$ and $l$ have the same parity. 
\label{rem_parity}
\end{rem}

We shall also need the following deep vanishing statement, due to J.-M. Bismut. This statement is 
what makes the calculations below possible and it would be be very interesting to have a better conceptual understanding of it. Its proof relies on the comparison between two completely different types of analytic torsion 
(holomorphic torsion and flat torsion) and it can be vaguely understood as a compatibility between the two sides of a kind of Hilbert correspondence. 
\begin{theor}
Let $h:M\to N$ be a proper and smooth morphism of complex manifolds. Let 
$g$ be a finite automorphism of $M$ over $N$ (ie $g$ acts fiberwise). Suppose 
that $h$ is endowed with a $g$-invariant Kähler fibration structure 
$\omega_h$. Then the element
$$
\sum_{k\geq 0}(-1)^k T_g(\mtr{\Omega}_h^k)\in \widetilde{\frak A}(N)
$$
vanishes.
\label{thBV}
\end{theor}
\beginProof
See \cite{Bismut-HdR}. 
\endProof
\begin{rem}
\rm In the paper \cite{BFL-Hol} it is shown that in the non-equivariant setting
the vanishing property stated in Theorem \ref{thBV} can be used to characterise 
the holomorphic torsion form axiomatically.\end{rem}

Let $f:X\to Y$ be a $\mn$-equivariant smooth and projective morphism 
of regular arithmetic varieties over an arithmetic ring $D$. Let $g$ be the automorphism 
of $X$  corresponding to $\zeta\in\mn(\mC)$. Suppose that $X(\mC)$ is endowed with a $g$-invariant 
Kähler fibration structure $\omega_f$ with respect to $f(\mC)$ and suppose that the $\mn$-structure of $Y$ is trivial. 
Suppose also that there is an equivariant line bundle on $X$, which is ample relatively to $f$.

We shall apply Theorem \ref{EAGRR} to the elements of the relative de Rham complex of $f$. To ease notation, let us write $\mtr{H}^k_\Db(X/Y)$ for the hermitian equivariant vector bundle
$$
\mtr{H}^k_\Db(X/Y):=\bigoplus_{p+q=k}\R^q f_*(\mtr{\Omega}_f^p)
$$
and $H^{p,q}(X/Y)$ for the vector bundle
$$
H^{p,q}(X/Y):=\R^q f_*({\Omega}_f^p).
$$

Theorem \ref{EAGRR} together with the Borel-Serre identity (see the end of section \ref{secERRT}) now gives the identity
\begin{equation}
\sum_k(-1)^k\ari{\ch}_\mn(\mtr{H}^k_\Db(X/Y))=f^{\mn}_*(c^\top(\mtr{\T f}))-\int_{X(\mC)_g/Y(\mC)}c^\top(\T f)R_g(\T f)
\label{afe}
\end{equation}
in $\ari{\CHOW}^\bullet_{\mQ(\mn)}(Y)$. 
In particular, for any $l\geq 1$, 
\medskip

\begin{equation}
\boxed{
\sum_k(-1)^k\ari{\ch}^{[l]}_\mn(\mtr{H}^k_\Db(X/Y))=-\int_{X(\mC)_g/Y(\mC)}c^\top(\T f)R_g^{[l-1]}(\T f)}\label{AGB}
\end{equation}

\medskip
We shall see below that equation \refeq{AGB} carries astonishingly deep arithmetic information. It should be viewed as a lifting to Arakelov theory of the relative equivariant form of the Gauss-Bonnet formula.

We shall now translate equation \refeq{AGB} into a statement about 
logarithmic derivatives of Dirichlet $L$-functions at negative integers. That this 
kind of translation should be possible is suggested by Lemma \ref{lemLChi} and the definition 
of the $R_g$-genus. 

We compute

\begin{lemma} For any $l\geq 1$ we have
$$
\sum_{\sigma\in\uZn}\wt{R}_{g^\sigma}^{[l-1]}(\T f)\chi(\sigma)=
\tau(\chi)\Big[2L'(\bar\chi,1-l)+\CH_{l-1}\cdot L(\bar\chi,1-l)\Big]\sum_{u\in\Zn}\ch^{[l-1]}(\T f_u)\bar\chi(u)
$$
\label{RgFourierLem}
\end{lemma}
\beginProof This is an immediate consequence of Lemma \ref{lemLChi} and the definition of the 
class $\wt{R}_g$.\endProof
\begin{lemma}
Let $h:M\to N$ be a projective and smooth morphism between quasi-projective complex manifolds. Let 
$g$ be a finite automorphism of $M$ over $N$ (ie $g$ acts fiberwise). Let 
\mbox{$h_0:M_g\to N$} be the induced morphism (which is smooth). The equality of 
characteristic classes in Aeppli cohomology
\begin{align*}
\sum_{u\in\Zn}\sum_{p,q}(-1)^{p+q}p\cdot \ch^{[l]}(&H^{p,q}(M/N)_u)\bar\chi(u)
\\
&=
-L(\bar\chi,-l)\int_{M_g/N}\sum_{u\in\Zn}c^\top(\T h_0)\ch^{[l]}((\Omega_h|_{M_g})_u)\bar\chi(u)
\end{align*}
holds.
\label{VMappLem}
\end{lemma}
Here $H^{p,q}(M/N)_u$ (resp. $(\Omega_h|_{M_g})_u$ is the largest subbundle of $H^{p,q}(M/N)$ (resp. $\Omega_h|_{M_g}$) where 
$g$ acts by multiplication by $\zeta^u$.
\beginProof This is an immediate consequence of Lemma \ref{lemLChi} and Corollary \ref{corGlob}.
\endProof

If $K\subseteq\mC$ is a subfield and $\chi$ is a Dirichlet character, we shall write 
$K(\chi)$ for the subfield of $\mC$ obtained by adjoining all the values of $\chi$ to 
$K$. 

Combining Lemma \ref{RgFourierLem} with 
equality \refeq{AGB}, we get the following.
For any primitive \mbox{Dirichlet} character 
modulo $n$, the equality: 

\begin{multline}
\sum_k(-1)^k\sum_{u\in\Zn}\ari{\ch}^{[l]}(\mtr{H}^k_\Db(X/Y)_u)\bar\chi(u)
\\
=-{1\over 2}\Big[2L'(\bar\chi,1-l)-2L'(\bar\chi,1-l)(-1)^{\chi(-1)+l-1}
+\CH_{l-1}\cdot \Big(L(\bar\chi,1-l)-L(\bar\chi,1-l)(-1)^{\chi(-1)+l-1}\Big)\Big]
\\
\cdot\int_{X(\mC)_g/Y(\mC)}c^\top(\T f)\sum_{u\in\Zn}\ch^{[l-1]}(\T f_u)\bar\chi(u)\,\,\,\,\,\,\,\,\,\,\,\label{eqAT}
\end{multline}
holds in $\ari{\CHOW}^l_{\mQ(\mu_{n})(\chi)}(Y)$.

Notice that if $\chi$ and $l$ have the same parity, we have
$$
L(\bar\chi,s)-L(\bar\chi,s)(-1)^{\chi(-1)+l-1}=2L(\bar\chi,s)
$$
whereas if $\chi$ and $l$ do not have the same parity then
$$
L(\bar\chi,s)-L(\bar\chi,s)(-1)^{\chi(-1)+l-1}=0.
$$
So we obtain from the equation \refeq{eqAT}: if $\chi$ and $l$ have the same parity, then
\begin{eqnarray*}
&&\sum_k(-1)^k\sum_{u\in\Zn}\ari{\ch}^{[l]}(\mtr{H}^k_\Db(X/Y)_u)\bar\chi(u)
\\&=&-\Big[2L'(\bar\chi,1-l)+\CH_{l-1}\cdot L(\bar\chi,1-l)\Big]\int_{X_g(\mC)/Y(\mC)}c^\top(\T f)\sum_{u\in\Zn}\ch^{[l-1]}(\T f_u)\bar\chi(u)\\
&=& -\Big[2L'(\bar\chi,1-l)+\CH_{l-1}\cdot L(\bar\chi,1-l)\Big]\int_{X_g(\mC)/Y(\mC)}c^\top(\T f)\sum_{u\in\Zn}\ch^{[l-1]}(\T f_{-u})\bar\chi(-u)\\
&=&
 -\Big[2L'(\bar\chi,1-l)+\CH_{l-1}\cdot L(\bar\chi,1-l)\Big]\int_{X_g(\mC)/Y(\mC)}c^\top(\T f)\sum_{u\in\Zn}\ch^{[l-1]}((\Omega_{f,u})^\vee)\bar\chi(-u)\\
 &=&
 -(-1)^{l-1}\bar\chi(-1)
 \Big[2L'(\bar\chi,1-l)+\CH_{l-1}\cdot L(\bar\chi,1-l)\Big]\\
 & &\hspace{6.5cm}\int_{X_g(\mC)/Y(\mC)}c^\top(\T f)\sum_{u\in\Zn}\ch^{[l-1]}(\Omega_{f,u})\bar\chi(u)\\
 &=&
 \Big[2L'(\bar\chi,1-l)+\CH_{l-1}\cdot L(\bar\chi,1-l)\Big]\int_{X_g(\mC)/Y(\mC)}c^\top(\T f)\sum_{u\in\Zn}\ch^{[l-1]}(\Omega_{f,u})\bar\chi(u)
\end{eqnarray*}
and if $\chi$ and $l$ do not have the same parity, then
\begin{equation}
\sum_k(-1)^k\sum_{u\in\Zn}\ari{\ch}^{[l]}(\mtr{H}^k_\Db(X/Y)_u)\bar\chi(u)=0.
\label{collapse}
\end{equation}
Finally, if $\chi$ and $l$ have the same parity (hence $L(\bar\chi,1-l)\not=0$ by Remark \ref{rem_parity}) then we obtain using 
Lemma \ref{VMappLem} that 
\begin{align}
\sum_k(-1)^k&\sum_{u\in\Zn}\ari{\ch}^{[l]}(\mtr{H}^k_\Db(X/Y)_u)\bar\chi(u)
\nonumber\\
=&-\sum_k(-1)^k\Big[2{L'(\bar\chi,1-l)\over L(\bar\chi,1-l)}+\CH_{l-1}\Big]\sum_{u\in\Zn}\sum_{p+q=k}p\cdot \ch^{[l-1]}(H^{p,q}(X/Y)_u)\bar\chi(u)
\label{AFUNDEQ}
\end{align}

Note that this equality does not depend on the initial choice of root of unit $\zeta_n$ anymore.

Now suppose that $\chi$ is not primitive. Then we apply \refeq{AFUNDEQ} again, but 
replace the action of $\mu_n$ by the action of its subgroup scheme 
$\mu_{f_\chi}$, where $f_\chi$ is the conductor of $\chi$. We shall write $\chi_{\rm Prim}$ 
for the primitive character modulo $f_\chi$ associated with $\chi$. 
Replacing $\chi$ by $\bar\chi$ for convenience, we finally get the following basic formula. 
We shall encase it in a theorem to underline its importance.

\begin{theor}
Let $f:X\ra Y$ be an $\mn$-equivariant smooth and projective morphism of equivariant regular arithmetic 
varieties. Suppose that the $\mn$-action on $Y$ is trivial. Fix a $\mn(\mC)$-invariant Kähler fibration structure $\omega_f$ for $f$ on $X$ and suppose that there is a $\mn$-equivariant line bundle on 
$X$, which is ample relatively to $f$. Let $\chi$ be a Dirichlet character modulo $n$. Then the equation
\begin{eqnarray}
&&\sum_k(-1)^k\sum_{u\in\Zn}\ari{\ch}^{[l]}(\mtr{H}^k_\Db(X/Y)_u)\chi_{\rm Prim}(u)
\nonumber\\&=&-\sum_k(-1)^k\Big[2\,{L'(\chi_{\rm Prim},1-l)\over L(\chi_{\rm Prim},1-l)}+\CH_{l-1}\Big]\sum_{u\in\Zn}\sum_{p+q=k}p\cdot \ch^{[l-1]}(H^{p,q}(X/Y)_u)\chi_{\rm Prim}(u)\,\,\,\,\,\,\,\,\,\,\label{AGBT}
\end{eqnarray}
holds in $\ari{\CHOW}_{\mQ(\mu_{f_\chi})(\chi_\Prim)}(Y)$, if $\chi$ and $l$ have the same parity (hence $L(\chi_{\rm Prim},1-l)\not=0$ by Remark \ref{rem_parity}). If $\chi$ and $l$ do not have the same parity then 
\begin{eqnarray*}
\sum_k(-1)^k\sum_{u\in\Zn}\ari{\ch}^{[l]}(\mtr{H}^k_\Db(X/Y)_u)\chi_{\rm Prim}(u)=0
\end{eqnarray*}
\label{AGBF}
\end{theor}
This is the promised translation of formula \refeq{AGB}. 

Now notice that in Theorem \ref{AGBF}, it is very natural to wonder whether the equality 
holds before the alternating sum  $\sum_k(-1)^k$ is taken on both sides. It is difficult 
to make a meaningful conjecture about this 'separation of weights' (in particular because 
the K\"ahler fibration structure $\omega_f$ is defined on $X$ and not only on the Gauss-Manin bundles). It nevertheless makes sense to conjecture the following purely geometric 
statement.

\begin{conj}
Let $f:X\ra Y$ be an $\mn$-equivariant smooth and projective morphism of equivariant regular arithmetic 
varieties. Suppose that the $\mn$-action on $Y$ is trivial. Let $\chi$ be a Dirichlet character modulo $n$. Then for any $k\geq 1$ the equation
\begin{eqnarray*}
\sum_{u\in\Zn}{\ch}^{[l]}({H}^k_\Db(X/Y)_u)\chi_{\rm Prim}(u)=0
\end{eqnarray*}
holds in ${\CHOW}^l(Y)_{\mQ(\mu_{f_\chi})(\chi_\Prim)}$ if $\chi$ and $l$ have the same parity. 
\label{AGBFgeom}
\end{conj}

When $D=\mC$ and $\chi=1$, this conjecture was studied and refined in \cite[Conj. 1.1]{MR-Conj}. See also \cite{EV-Chern} for this conjecture. 

In the direction of 'separation of weights' in the context of Arakelov geometry, we can nevertheless prove prove the following result.

\begin{prop}
Let $f:\CA=X\ra Y$ be an $\mn$-equivariant smooth and projective morphism of equivariant regular arithmetic 
varieties and suppose that $\CA$ is an abelian scheme. Suppose that there is a $\mn$-equivariant line bundle on 
$X$, which is ample relatively to $f$. Suppose that the $\mn$-action on $Y$ is trivial and that $\CA_\mn$ is finite over $Y$. Fix a $\mn(\mC)$-invariant Kähler fibration structure $\omega_f$ and suppose also that $\omega_f$ is translation invariant on the fibres of $f(\mC)$ and that $${1\over\dim(\CA/Y)!}\int_{\CA(\mC)/Y(\mC)}\omega_f^{\dim(\CA/Y)}=1.$$ Let $\chi$ be a Dirichlet character modulo $n$. Then the equation
\begin{eqnarray}
&&\sum_{u\in\Zn}\ari{\ch}^{[l]}(\mtr{H}^k_\Db(X/Y)_u)\chi_{\rm Prim}(u)
\nonumber\\&=&-\Big[2{L'(\chi_{\rm Prim},1-l)\over L(\chi_{\rm Prim},1-l)}+\CH_{l-1}\Big]\sum_{u\in\Zn}\sum_{p+q=k}p\cdot \ch^{[l-1]}(H^{p,q}(X/Y)_u)\chi_{\rm Prim}(u)\,\,\,\,\,\,\,\,\,\,
\end{eqnarray}
holds in $\ari{\CHOW}^l_{\mQ(\mu_{f_\chi})(\chi_\Prim)}(Y)$ if $\chi$ and $l$ have the same parity (hence $L(\chi_{\rm Prim},1-l)\not=0$ by Remark \ref{rem_parity}). 
\label{propCTA}
\end{prop}
To prove this, we shall need the following combinatorial lemmata.

Consider the following formal power series:
\[
\exp(x) := \sum_{j=0}^\infty \frac{x^{j}}{j!}
\]
and
\[
\log(1 + x):=\sum_{j=1}^\infty (-1)^{j+1}\frac{x^j}{j}.
\]
 
\begin{lemma}
Let $\lambda\in\mn(\mC)$, $\lambda\not=1$. Then the equality
$$
\log({1-\lambda\cdot\exp(x)\over 1-\lambda})=
 -\sum_{j\geq 1}\zeta_L(\lambda,1-j){x^j\over j!}
$$
holds in $\mC[[x]]$. 
\label{MainComb}
\end{lemma}
\beginProof See \cite[Lemma 4]{MR-Order}. 
\endProof

\begin{lemma}
Let $\mtr{V}$ be a $\mn$-equivariant hermitian bundle on 
an arithmetic variety $Z$. Suppose that the $\mn$-action on $Z$ is trivial.  
Then we have the equality
\begin{equation}
\log\Big[\Big(\prod_{u=0}^{n-1}(1-\zeta^u)^{\rk(V_u)}\Big)^{-1}\sum_{r\geq 0}(-1)^r\ach_\mn(\Lambda^r(\mtr{V}))\Big]=
-\sum_{l\geq 1}\sum_{u=0}^{n-1}\zeta_L(\zeta^u,1-l)\ach^{[l]}(\mtr{V}_u)
\label{eqMC2}
\end{equation}
in $\ari{\CHOW}_{\mQ(\mu_n)}(Z)$.
\label{lemMC2}
\end{lemma}
\beginProof Notice first that if $\mtr{V}$ and $\mtr{W}$ 
are $\mn$-equivariant vector bundles, then we have 
$$
\sum_{r\geq 0}(-1)^r\ach_\mn(\Lambda^r(\mtr{V}\oplus\mtr{W}))=
\sum_{r\geq 0}(-1)^r\ach_\mn(\Lambda^r(\mtr{V}))\cdot 
\sum_{r\geq 0}(-1)^r\ach_\mn(\Lambda^r(\mtr{W})).
$$
Thus
\begin{eqnarray*}
&&\log\Big[\Big(\prod_{u=0}^{n-1}(1-\zeta^u)^{\rk((V\oplus W)_u)}\Big)^{-1}\sum_{r\geq 0}(-1)^r\ach_\mn(\Lambda^r(\mtr{V}\oplus\mtr{W}))\Big]\\&=&
\log\Big[\Big(\prod_{u=0}^{n-1}(1-\zeta^u)^{\rk(V_u)}\Big)^{-1}\sum_{r\geq 0}(-1)^r\ach_\mn(\Lambda^r(\mtr{V}))\Big]\\&~&\hspace{2cm}+\log\Big[\Big(\prod_{u=0}^{n-1}(1-\zeta^u)^{\rk(W_u)}\Big)^{-1}\sum_{r\geq 0}(-1)^r\ach_\mn(\Lambda^r(\mtr{W}))\Big]
\end{eqnarray*}
and in particular both sides of the equality \refeq{eqMC2} are additive in $\mtr{V}$. 
By the splitting principle, we are thus reduced to the case of a line bundle, in which case 
the lemma reduces to Lemma \ref{MainComb}.
 \endProof

\beginProof ({Proof of Proposition \ref{propCTA}})  Notice first that from the definition of the $L_2$-metric and the assumption that 
${1\over\dim(\CA/Y)!}\int_{\CA(\mC)/Y(\mC)}\omega_f^{\dim(\CA/Y)}=1$,  
there exists an isometric isomorphism
$$
\Lambda^r(\mtr{H}^1_\Db(\CA/Y))\simeq \mtr{H}^r_\Db(\CA/Y)
$$
for all $r\geq 0$. We now apply equality \refeq{afe} to
$X=\CA$ over $Y$. We obtain
$$
\sum_{r\geq 0}(-1)^r\ach_\mn(\Lambda^r(\mtr{H}^1_\Db(\CA/Y)))=
\Big(\prod_{u=0}^{n-1}(1-\zeta^u)^{\rk({H}^1(\CA/Y)_u)}\Big)(1-R_g(f_*(\T f)))
$$
Applying Lemma \ref{lemMC2}, we obtain
$$
\sum_{l\geq 1}\sum_{u=0}^{n-1}\zeta_L(\zeta^u,1-l)\ach^{[l]}(\mtr{H}^1_\Db(\CA/Y)_u)
=R_g(f_*(\T f))
$$
Now there is an equivariant isomorphism $f_*(\T f(\mC))\simeq H^{1,0}(X/Y)(\mC)^\vee$ 
(given by the polarisation induced by a $\mn$-equivariant relatively ample line bundle) and applying Lemma \ref{lemLChi}, we finally obtain that
\begin{eqnarray*}
&&\sum_{l\geq 1}\sum_{u=0}^{n-1}L(\bar\chi,1-l)\ach^{[l]}(\mtr{H}^1_\Db(\CA/Y)_u)\bar\chi(u)\\&=&
(-1)^{l-1}\bar\chi(-1)\sum_{l\geq 1}{1\over 2}\Big[2L'(\bar\chi,1-l)-2L'(\bar\chi,1-l)(-1)^{\sgn(\chi)+l-1}\\&+&\CH_{l-1}\cdot \Big(L(\bar\chi,1-l)-L(\bar\chi,1-l)(-1)^{\sgn(\chi)+l-1}\Big)\Big]
\cdot\ch^{[l-1]}(H^{1,0}(X/Y))\bar\chi(u)
\end{eqnarray*}
for any primitive Dirichlet character $\chi$ modulo $n$. We can now conclude 
following the same line of argument as in the proof of Theorem \ref{AGBF}. \endProof

\begin{cor}
Conjecture \ref{AGBFgeom} holds if $X$ is an abelian scheme over $Y$, $k=1$ and $X_\mn$ is finite over $Y$.
\end{cor}

We now wish to translate Proposition \ref{propCTA} into the language of complex 
multiplication of abelian schemes. For this and later applications, we shall need the following 

\begin{lemma}
Suppose that $L,K$ are number fields and that all the embeddings of 
$K$ into $\mC$ factor through an embedding of $L$ into $\mC$. Let 
$d_K$ be the discriminant of $K$. Then there is a canonical isomorphism 
of $\CO_L$-algebras
$$
(\CO_L\otimes\CO_K)[{1\over d_K}]\simeq 
\bigoplus_{\sigma:K\hookrightarrow L}\CO_L[{1\over d_K}]
$$
such that $l\otimes k\mapsto\oplus_\sigma\, l\cdot \sigma(k).$
\label{dclem}
\end{lemma}
\beginProof
Notice to begin with that we have an isomorphism of $L$-algebras
\begin{equation}
L\otimes_\mQ K\simeq \bigoplus_{\sigma:K\hookrightarrow L} 
L
\label{geniso}
\end{equation}
such that $l\otimes k\mapsto\oplus_\sigma\, l\cdot \sigma(k).$
This can be seen by writing $K\simeq\mQ[t]/(P(t))$ for some 
monic irreducible polynomial $P(t)$ and noticing that by assumption, $P(t)$ splits in $L$. 

Notice now that $\Spec\,\CO_K[{1\over d_K}]\to\Spec\ \mZ[{1\over d_K}]$ is by construction 
a finite and \'etale morphism. Hence the morphism 
$$
\Spec\,\CO_L\otimes\CO_K[{1\over d_K}]\to\Spec\,\CO_L[{1\over d_K}]
$$
is also finite and \'etale. Thus $\Spec\,\CO_L\otimes\CO_K[{1\over d_K}]$ is the disjoint 
union of its irreducible components and any of these components, say $C$, is integral, 
finite and \'etale over $\Spec\,\CO_L[{1\over d_K}]$. 
On the other hand, notice that the morphism $C_L\to\Spec\,L$ is an isomorphism 
because of the existence of the decomposition \refeq{geniso}. Thus there is a section $\Spec\,L\to C_L$, which 
extends uniquely to a section $$\Spec\,\CO_L[{1\over d_K}]\to C$$ by the valuative 
criterion of properness and this section is an open immersion because 
$C\to\Spec\,\CO_L[{1\over d_K}]$ is \'etale. Hence this section is an isomorphism, since 
$C$ is integral. To summarise, the irreducible components of $\Spec\,\CO_L\otimes\CO_K[{1\over d_K}]$ are all images of sections of the morphism $
\Spec\,\CO_L\otimes\CO_K[{1\over d_K}]\to\Spec\,\CO_L[{1\over d_K}].
$
Furthermore, every section of the morphism 
$
\Spec\,L\otimes_\mQ K\to\Spec\,L
$
extends uniquely to a section over $\Spec\,\CO_L[{1\over d_K}]$, which is an open immersion and 
whose image is an irreducible component. 
Translating these two statements back into the language of rings gives the lemma.\endProof

Suppose now that $X=\CA$ is an abelian scheme over $Y$ and that there is an injection $\CO_{\mQ(\mn)}\hookrightarrow\End_Y(\CA)$ for some $n>1$. 
Suppose also the $n$ is invertible in the arithmetic base ring $D$, that there is a primitive $n$-th root of unity in $D$ and that $D$ is a localisation of the rings of integers of a number field. Then 
$\mu_{n,D}$ is isomorphic to the constant group scheme over $D$ associated with $\mZ/(n)$. We fix such an isomorphism; this is equivalent to choosing a primitive 
root of unity in $D$, or in other words to choosing an embedding 
$\iota:\CO_{\mQ(\mn)}\hookrightarrow D$. We are now given a $\mn$-action on $\CA$ over $Y$. Note that $n=2$ is allowed; in that case the $\mu_2$-action given by the injection $\CO_{\mQ(\mu_2)}=\mQ\hookrightarrow\End_Y(\CA)$ is given by the action of 
the automorphism $[-1]_\CA$. By Lemma 
\ref{dclem}, we have 
$\mtr{H}^1_\Db(\CA/Y)_u=0$ if $u$ is not prime to $n$. In particular, 
$\CA_\mn$ is finite over $Y$. Hence, for any Dirichlet character modulo $n$, 
Proposition \ref{propCTA} gives us the  
equality:
\begin{align*}
\sum_{u\in(\Zn)^*}\ari{\ch}^{[l]}(\mtr{H}^1_\Db(\CA/&Y)_u)\chi(u)
\\
=&-\Big[2\,{L'(\chi_{\rm Prim},1-l)\over L(\chi_{\rm Prim},1-l)}+\CH_{l-1}\Big]\sum_{u\in(\Zn)^*}\ch^{[l-1]}(H^{1,0}(\CA/Y)_u)\chi(u)\\
\end{align*}
in $\ari{\CHOW}^l_{\mQ(\mu_{f_\chi})(\chi_\Prim)}(Y)$ if $\chi$ and $l$ have the same parity (hence $L(\chi_{\rm Prim},1-l)\not=0$, see Remark \ref{rem_parity}.) 
This can be rewritten as:
\begin{align}
\sum_{\tau\in\Gal(\mQ(\mn)|\mQ)}&\ari{\ch}^{[l]}(\mtr{H}^1_\Db(\CA/Y)_{\iota\circ\tau})\chi(\tau)
\nonumber
\\
=&-\Big[2\,{L'(\chi_{\rm Prim},1-l)\over L(\chi_{\rm Prim},1-l)}+\CH_{l-1}\Big]\sum_{\tau\in\Gal(\mQ(\mn)|\mQ)}\ch^{[l-1]}(H^{1,0}(\CA/Y)_{\iota\circ\tau})\chi(\tau)\label{BRR}
\end{align}
where $\mtr{H}^1_\Db(\CA/Y)_{\iota\circ\tau}$ is the subsheaf on which $\CO_{\mQ(\mn)}$ acts 
via $\iota\circ\tau$. Notice that since $n$ is invertible in $D$, Lemma \ref{dclem} implies that there is an inner direct sum
$$
\bigoplus_{\tau\in\Gal(\mQ(\mn)|\mQ)}{H}^1_\Db(\CA/Y)_{\iota\circ\tau}\simeq 
{H}^1_\Db(\CA/Y).
$$
\begin{rem}\rm Notice the interesting fact that the truth of equation \refeq{BRR} 
is independent of the embedding $\iota:\CO_{\mQ(\mn)}\hookrightarrow D$. Indeed 
if $\iota_1:\CO_{\mQ(\mn)}\hookrightarrow D$ is another embedding then there exist 
$\tau_1\in \Gal(\mQ(\mn)|\mQ)$ such that $\iota_1=\iota\circ\tau_1$ (because $\mQ(\mn)$ 
is a Galois extension of $\mQ$). Thus 
\begin{multline*}
\sum_{\tau\in\Gal(\mQ(\mn)|\mQ)}\ari{\ch}^{[l]}(\mtr{H}^1_\Db(\CA/Y)_{\iota_1\circ\tau})\chi(\tau)
\\
=\sum_{\tau\in\Gal(\mQ(\mn)|\mQ)}\ari{\ch}^{[l]}(\mtr{H}^1_\Db(\CA/Y)_{\iota\circ\tau})\chi(\tau_1^{-1}\circ\tau)\\
=
\chi(\tau_1^{-1})\left(\sum_{\tau\in\Gal(\mQ(\mn)|\mQ)}\ari{\ch}^{[l]}(\mtr{H}^1_\Db(\CA/Y)_{\iota\circ\tau})\chi(\tau)\right)
\end{multline*}
and similarly
\begin{eqnarray*}
&&-\Big[2\,{L'(\chi_{\rm Prim},1-l)\over L(\chi_{\rm Prim},1-l)}+\CH_{l-1}\Big]\sum_{\tau\in\Gal(\mQ(\mn)|\mQ)}\ch^{[l-1]}(H^{1,0}(\CA/Y)_{\iota_1\circ\tau})\chi(\tau)\\
&=&
-\Big[2\,{L'(\chi_{\rm Prim},1-l)\over L(\chi_{\rm Prim},1-l)}+\CH_{l-1}\Big]\sum_{\tau\in\Gal(\mQ(\mn)|\mQ)}\ch^{[l-1]}(H^{1,0}(\CA/Y)_{\iota\circ\tau})\chi(\tau_1^{-1}\circ\tau)\\
&=&
-\chi(\tau_1^{-1})\left(\Big[2\,{L'(\chi_{\rm Prim},1-l)\over L(\chi_{\rm Prim},1-l)}+\CH_{l-1}\Big]\sum_{\tau\in\Gal(\mQ(\mn)|\mQ)}\ch^{[l-1]}(H^{1,0}(\CA/Y)_{\iota\circ\tau})\chi(\tau)\right)
\end{eqnarray*}
and we can thus conclude that if equality \refeq{BRR} is true for 
a certain embedding $\iota$ then it is true for any such embedding. This might seem a moot point since we know anyway that equality \refeq{BRR} is true but it seemed 
worth underlining in view of the following.
\end{rem}

Equality \refeq{BRR} suggests the following conjecture:

\begin{conj}
Suppose that $K$ is a finite Galois extension of $\mQ$. Suppose that there is an element $c\in\Gal(K|\mQ)$ in the center of $\Gal(K|\mQ)$ such that for all embeddings $\iota:K\to\mC$ and all $k\in K$, we have $\iota(c(k))=\mtr{\iota(k)}$ (where $\mtr{(\cdot)}$ refers to complex conjugation). Suppose that all the embeddings of $K$ into $\mC$ factor through an embedding of ${\rm Frac}(D)$ into $\mC$. Suppose finally that the discriminant of $K$ is invertible in $D$ and that 
$D$ is a localisation of the ring of integers of a number field. 

Let $f:\CA\to Y$ be an abelian scheme and suppose that 
we are given an embedding of rings $\rho:\CO_K\hookrightarrow\End_Y(\CA).$ 

Let 
$\chi:\Gal(K|\mQ)\to\mC$ be an irreducible Artin character and let $l\geq 1$. 
 
Suppose given a K\"ahler fibration structure $\nu_f$  such that 

$\bullet$ $\nu_f$ represents the first Chern class of a relatively ample line bundle;

$\bullet$ for any $x\in\CO_K$, the endomorphism $\rho(x)^*$ of 
${H}^1_\Db(\CA/Y)(\mC)$ is adjoint to the endomorphism $\rho(c(x))^*$ of ${H}^1_\Db(\CA/Y)(\mC)$, with respect 
to the metric coming from $\nu_f$. 

Suppose that $\chi$ and $l$ have the same parity (hence $L(\chi,1-l)\not=0$ by Remark \ref{rem_parity}).

Then for any embedding $\iota:\CO_K\hookrightarrow D$ we have:
\begin{align*}
\sum_{\tau\in\Gal(K|\mQ)}\ari{\ch}^{[l]}(\mtr{H}^1_\Db(\CA/&Y)_{\iota\circ\tau})\chi(\tau)
\nonumber\\ =&-\Big[2\,{L'(\chi,1-l)\over L(\chi,1-l)}+\CH_{l-1}\Big]\sum_{\tau\in\Gal(K|\mQ)}\ch^{[l-1]}(H^{1,0}(\CA/Y)_{\iota\circ\tau})\chi(\tau)\nonumber
\end{align*}
in $\ari{\CHOW}^l_{\bar\mQ}(Y).$ 
\label{RIconj}
 \end{conj}
 The endomorphism $\rho(x)^*$ is the endomorphism of ${H}^1_\Db(\CA/Y)(\mC)$ obtained by pull-back. The sheaves $\mtr{H}^1_\Db(\CA/Y)$ are understood to carry the $L_2$-metric induced by the K\"ahler fibration structure $\nu_f$.
The notation $L(\chi,s)$ refers to the Artin $L$-function associated with $\chi$ (see 
\cite{Tate-Stark} or \cite[\S\ 7.10-12]{Neu} for an introduction). 
 Note that if $K=\mQ(\mn)$ then $\chi$ can be identified with a Dirichlet character $\chi_0$
 via the canonical isomorphism $\Gal(\mQ(\mn)|\mQ)\simeq(\mZ/(n))^*$ and then one 
 has $L(\chi,s)=L(\chi_{0,\Prim},s)$. 
 
\begin{rem}\rm If $Y=\Spec\,D$, $K$ is a CM field and the generic fibre $f:\CA\to Y$ is an abelian variety of dimension ${1\over 2}[K:\mQ]$ (in particular  
the generic fibre of $\CA$ has CM by $\CO_K$), then there always exists a 
K\"ahler fibration structure of the type described in Conjecture \ref{RIconj}. 
See \cite[after Th. A]{Rosen-Ab}.\end{rem}
\begin{rem} \rm If $K=\mQ(\mn)$ then a polarisation with the 
properties required in Conjecture \ref{RIconj} can be constructed as follows. 
Choose first a $\mn(\mC)$-equivariant relatively ample line bundle $\CL$ on $\CA(\mC).$ 
Such a line bundle can be obtained in the following way. Let 
$\CM$ be a relatively ample line bundle on $\CA(\mC)$ (without equivariant structure). 
The line bundle $\otimes_{a\in\mn(\mC)}a^*\CM$ then carries a $\mn(\mC)$-equivariant 
structure and is also relatively ample. Suppose without restriction of generality that the restriction of $\CL$ to the $0$-section of $\CA(\mC)$ is 
an equivariantly trivial line bundle and choose a trivialisation. There is then a unique hermitian 
metric on $\CL$, whose first Chern character form is translation invariant on the fibres 
of $\CA(\mC)\to Y(\mC)$ and such that the trivialising map has norm $1$ 
(see eg \cite[II, 2.1]{MB} for all this). Let $\mtr{\CL}$ be the resulting hermitian line bundle. 
The first Chern character form $c^1(\mtr{\CL})$ is then a $\mn(\mC)$-invariant Kähler fibration 
structure for $\CA(\mC)\to Y(\mC)$ and it satisfies the properties required 
in Conjecture \ref{RIconj} because for any $\zeta\in\mn(\mC)$ the adjoint of $\rho(\zeta)^*$ is 
then $\rho(\zeta)^{*,-1}=\rho(\zeta^{-1})^*=\rho(\bar\zeta)$ and $\mn(\mC)$ generates 
$\CO_{\mQ(\mn)}$ as a $\mZ$-module. 
\label{remcyc}
\end{rem}
\begin{rem}\rm 
Notice that the assumptions of Conjecture \ref{RIconj} imply that the subbundles 
$\mtr{H}^1_\Db(\CA/Y)_{\iota\circ\tau}$ of $\mtr{H}^1_\Db(\CA/Y)$ are orthogonal to each other. This follows from the fact that by construction the pull-back endomorphisms $\rho(x)^*$ commute 
with their adjoints for any $x\in\CO_K$.
\label{remorth}
\end{rem}

One might wonder how much Conjecture \ref{RIconj} depends on the polarisation.
We shall show
\begin{prop}
Suppose that $K$ is a finite Galois extension of $\mQ$. Suppose that there is an element $c\in\Gal(K|\mQ)$ such that for all embeddings $\iota:K\to\mC$ and all $k\in K$, we have $\iota(c(k))=\mtr{\iota(k)}$ (where $\mtr{(\cdot)}$ refers to complex conjugation). Suppose that all the embeddings of $K$ into $\mC$ factor through an embedding of ${\rm Frac}(D)$ into $\mC$. Suppose finally that the discriminant of $K$ is invertible in $D$ and that 
$D$ is a localisation of the ring of integers of a number field. 

Let $f:\CA\to Y$ be an abelian scheme and suppose that 
we are given an embedding of rings $\rho:\CO_K\hookrightarrow\End_Y(\CA).$ 

Suppose given a K\"ahler fibration structure $\nu_f$ (resp. $\kappa_f$) such that 

$\bullet$ $\nu_f$ (resp. $\kappa_f$)  represents the first Chern class of a relatively ample line bundle;

$\bullet$ for any $x\in\CO_K$, the pull-back endomorphism $\rho(x)^*$ of 
${H}^1_\Db(\CA/Y)_\mC$ is adjoint to the pull-back endomorphism $\rho(c(x))^*$ of ${H}^1_\Db(\CA/Y)_\mC$, with respect 
to the metric coming from $\nu_f$ (resp. $\kappa_f$). 

Suppose that $l$ and $\chi$ have the same parity. Then for any embedding $\iota:\CO_K\hookrightarrow D$, we have
$$
\sum_{\tau\in\Gal(K|\mQ)}\ari{\ch}^{[l]}(\mtr{H}^1_\Db(\CA/Y)_{\nu_f,\iota\circ\tau})\chi(\tau)=
\sum_{\tau\in\Gal(K|\mQ)}\ari{\ch}^{[l]}(\mtr{H}^1_\Db(\CA/Y)_{\kappa_f,\iota\circ\tau})\chi(\tau)
$$
in $\ari{\CHOW}^l_{\bar\mQ}(Y)$.
\label{RIconjindep}
 \end{prop}
Here we write $\mtr{H}^1_\Db(\CA/Y)_{\nu_f,\iota\circ\tau}$ (resp. 
$\mtr{H}^1_\Db(\CA/Y)_{\kappa_f,\iota\circ\tau}$) for the bundle ${H}^1_\Db(\CA/Y)_{\iota\circ\tau}$ 
endowed with the $L_2$-metric induced by $\nu_f$ (resp. $\kappa_f$).

Proposition \ref{RIconjindep} in particular says that the truth of Conjecture \ref{RIconj} does not depend 
on the choice of the polarisation. To prove Proposition \ref{RIconjindep}, we shall 
need the following lemma.

\begin{lemma}
\label{lem724}
Let $M$ be a complex manifold and let $(V,h_0)$ be a holomorphic vector bundle 
$V$ on $M$, endowed with a hermitian metric $h_0$. Let $\phi:V\to V$ be an 
automorphism of vector bundles and suppose that $\phi$ is positive definite with respect to $h_0$ (on each fibre of $V$). Let $h_1$ be the hermitian metric on $V$ defined by the formula 
$h_1(v,w):=h(\phi(v),w)$ for any elements $v,w\in V$ which lie in the same fibre. 

Let 
$\wt{\ch}(V,h_1,h_2)\in\wt{\mathfrak A}(M)$ be the Bott-Chern secondary class of the 
exact sequence
$$
\bar{\cal E}:0\to V\stackrel{\rm Id}{\to}V\to 0
$$
where the first non zero term from the right carries the metric $h_1$ and the 
second non zero term from the right carries the metric $h_2$. 

Then the eigenvalues 
of $\phi$ are locally constant on $M$ and we have
$$
\wt{\ch}(V,h_1,h_2)=\sum_{t\in\mR_{>0}}\log(t)\ch((V_t,h_1|_{V_t}))
$$
where $V_t$ is the kernel of $\phi-t\cdot\Id$.
\end{lemma}
\beginProof
Since $\phi$ is self adjoint on the fibres of $V$ with respect to $h_1$, the coefficients 
of the polynomial $\det(\Id-x\cdot \phi)\in\CO(M)[x]$ are real valued holomorphic 
functions and they are thus locally constant. Furthermore, the automorphism 
$\phi$ is diagonalisable on each fibre of $V$ and thus we have a decomposition of $V$ 
$$
V\simeq \bigoplus_{t\in\mR_{>0}}V_t
$$
as an inner orthogonal direct sum of vector bundles. Furthermore, we have by construction
$$
h_2(v,w)=h_1(t\cdot v,w)
$$
for any elements of $V_t$ that lie in the same fibre. Hence we have
$$
\wt{\ch}(V,h_1,h_2)=\sum_{t\in\mR_{>0}}\wt{\ch}(V_t,h_1|_{V_t},t\cdot h_1|_{V_t}).
$$
Now we have 
$$
\wt{\ch}(V_t,h_1|_{V_t},t\cdot h_1|_{V_t})=\log(t)\cdot\ch((V_t,h_1|_{V_t})).
$$
See eg \cite[ex. on p. 22]{Faltings-Lectures} for this.
\endProof

\beginProof (of Proposition \ref{RIconjindep}) 
We start with some preliminary considerations. 
Let $M$ be a projective complex manifold and let $L$ be an ample line bundle on $M$. Let $\omega\in H^{2}(M,\mC)$ be the first Chern class of $L$ in complex Betti cohomology. Let 
$v,w\in H^1(M,\mC)$. We shall write $\bar{v}$ for the complex conjugate of 
$v$ and $v_{0,1}$ (resp. $v_{1,0}$) for the Hodge components of $v$ (and similarly for $w$). By the discussion preceding Lemma \ref{HDRLem}, we have the formula
$$
\langle v,w\rangle_{{\rm Hodge},L}=\int_M v\wedge *\bar w=
{i\over (\dim(M)-1)!}\int_M v\wedge \omega^{\dim(M)-1}\wedge(\bar{w}_{0,1}-\bar{w}_{1,0})
$$
for the Hodge metric on $H^1(M,\mC)$. Now choose another 
ample line bundle $J$, with first Chern class $\eta\in H^{2}(M,\mC)$ say. 
The maps 
$$
\bullet\wedge\omega^{\dim(M)-1}:H^1(M,\mC)\to H^{2\dim(M)-1}(M,\mC)
$$
and 
$$
\bullet\wedge\eta^{\dim(M)-1}:H^1(M,\mC)\to H^{2\dim(M)-1}(M,\mC)
$$
are both isomorphisms by the Hard Lefschetz theorem for singular cohomology. These isomorphisms also respect the underlying $\mQ$-rational Hodge structures. Hence there is a unique isomorphism of $\mQ$-rational Hodge structures 
$$
M=M(L,J):H^1(M,\mQ)\to H^1(M,\mQ)
$$
such that $$\langle M(v),w\rangle_{{\rm Hodge},L}=\langle v,w\rangle_{{\rm Hodge},J}.$$
for all $v,w\in H^1(M,\mC)$. 
Since both $\langle \bullet,\bullet\rangle_{{\rm Hodge},L}$ and 
$\langle\bullet,\bullet\rangle_{{\rm Hodge},J}$ are Hermitian metrics, the 
isomorphism $M$ is necessarily positive definite for the metric $\langle \bullet,\bullet\rangle_{{\rm Hodge},L}$. Now suppose furthermore that 
we are given endomorphisms $e,d:H^1(M,\mC)\to H^1(M,\mC)$ of $\mC$-vector spaces and suppose that $e$ and $d$ commute and that $d$ is the adjoint of $e$ with respect to $\langle \bullet,\bullet\rangle_{{\rm Hodge},L}$ \underline{and} with respect to $\langle \bullet,\bullet\rangle_{{\rm Hodge},J}$. Then we contend that $e$ commutes 
with $M$. Indeed from the assumptions on $d$ and $e$ we may compute
$$
\langle e(M(v)),w\rangle_{{\rm Hodge},L}=\langle M(v),d(w)\rangle_{{\rm Hodge},L}
$$
$$
\langle M(e(v)),w\rangle_{{\rm Hodge},L}=\langle M(v),d(w)\rangle_{{\rm Hodge},L}
$$
and since $v,w$ are arbitrary we conclude that $e\circ M=M\circ e$. 

Now let us return to the matter at hand. A straightforward generalisation of the preceding calculation to 
a relative setting shows that there an automorphism of vector bundles 
$$
M=M(\nu_f,\kappa_f):{H}^1_\Db(\CA/Y)(\mC)\to {H}^1_\Db(\CA/Y)(\mC)
$$
which is self adjoint with respect to the $L_2$-metric induced by $\nu_f$ and 
such that 
$$
\langle M(\bullet),\bullet\rangle_{L_2,\nu_f}=\langle \bullet,\bullet\rangle_{L_2,\kappa_f}.
$$
Furthermore, for any $x\in\CO_K$, in view of the assumptions on 
$\rho(x)^*$, we see that $\rho(x)^*$ commutes with $M$. Thus 
$M$ respects the decomposition
$$
{H}^1_\Db(\CA/Y)(\mC)\simeq
\bigoplus_{\tau\in\Gal(K|\mQ)}{H}^1_\Db(\CA/Y)(\mC)_{\iota\circ\tau}.
$$
Thus, using Lemma \ref{lem724} we may compute
$$
\ari{\ch}^{[l]}(\mtr{H}^1_\Db(\CA/Y)_{\nu_f,\iota\circ\tau})-
\ari{\ch}^{[l]}(\mtr{H}^1_\Db(\CA/Y)_{\kappa_f,\iota\circ\tau})=
\sum_{t\in\mR_{>0}}\log(t)\ch^{[l-1]}(\mtr{H}^1_\Db(\CA/Y)(\mC)_{\iota\circ\tau,t})
$$
where ${H}^1_\Db(\CA/Y)(\mC)_{\iota\circ\tau,t}$ is the sub bundle 
of ${H}^1_\Db(\CA/Y)(\mC)_{\iota\circ\tau}$ corresponding to the eigenvalue $t$ of 
$M$. Now notice that ${H}^1_\Db(\CA/Y)(\mC)_{\iota\circ\tau,t}$ is isomorphic 
as a $C^\infty$-vector bundle to a flat bundle via the comparison isomorphism 
with the corresponding relative Betti cohomology sheaves. Hence 
$\ch(\mtr{H}^1_\Db(\CA/Y)(\mC)_{\iota\circ\tau,t})$ is $d$-exact in positive degrees 
and in particular the positive degree part of the expression
$$
\sum_{t\in\mR_{>0}}\log(t)\ch^{[l-1]}(\mtr{H}^1_\Db(\CA/Y)(\mC)_{\iota\circ\tau,t})
$$
vanishes in Aeppli cohomology. We conclude that the difference
$$
\sum_{\tau\in\Gal(K|\mQ)}\ari{\ch}^{[l]}(\mtr{H}^1_\Db(\CA/Y)_{\nu_f,\iota\circ\tau})\chi(\tau)-
\sum_{\tau\in\Gal(K|\mQ)}\ari{\ch}^{[l]}(\mtr{H}^1_\Db(\CA/Y)_{\kappa_f,\iota\circ\tau})\chi(\tau)
$$
vanishes if $l>1$. This settles the proposition for $l>1$. If $l=1$, the difference is
\begin{eqnarray*}
&&\sum_{\tau\in\Gal(K|\mQ)}\ari{\ch}^{1}(\mtr{H}^1_\Db(\CA/Y)_{\nu_f,\iota\circ\tau})\chi(\tau)-
\sum_{\tau\in\Gal(K|\mQ)}\ari{\ch}^{1}(\mtr{H}^1_\Db(\CA/Y)_{\kappa_f,\iota\circ\tau})\chi(\tau)\\
&=& 
\sum_{\tau\in\Gal(K|\mQ)}\sum_{t\in\mR_{>0}}\log(t)\rk({H}^1_\Db(\CA/Y)(\mC)_{\iota\circ\tau,t}))\chi(\tau)\\
&=&
\sum_{t\in\mR_{>0}}\log(t)\sum_{\tau\in\Gal(K|\mQ)}\rk({H}^1_\Db(\CA/Y)(\mC)_{\iota\circ\tau,t}))\chi(\tau).
\end{eqnarray*}
We shall now show that $\sum_{\tau\in\Gal(K|\mQ)}\rk({H}^1_\Db(\CA/Y)(\mC)_{\iota\circ\tau,t}))\chi(\tau)=0$ 
if $\chi$ is an odd character. This will conclude the proof of the proposition. 
To show this, we may suppose that $Y(\mC)$ is a finite set of points, so 
suppose for simplicity that $Y=\Spec\,\mZ$. In that case, 
${H}^1_\Db(\CA/Y)(\mC)_{\iota\circ\tau}$ reduces to a  
complex vector space. 
Via the comparison isomorphism, this vector space has a $\mQ$-rational 
structure and the automorphism $M$ respects this structure. 
Thus 
$$
\rk({H}^1_\Db(\CA/Y)(\mC)_{\iota\circ\tau,t}))=\rk({H}^1_\Db(\CA/Y)(\mC)_{\iota\circ c\circ\tau,t}))
$$
so that the function $\rk({H}^1_\Db(\CA/Y)(\mC)_{\iota\circ\tau,t}))$ is an even function 
on $\Gal(K|\mQ)$. We conclude that 
$$
\sum_{\tau\in\Gal(K|\mQ)}\rk({H}^1_\Db(\CA/Y)(\mC)_{\iota\circ\tau,t}))\chi(\tau)=0.
$$\endProof

We shall now prove
\begin{theor} Suppose that the assumptions of Conjecture \ref{RIconj} hold. 
Suppose furthermore that $K$ is an abelian extension of $\mQ$. Then the conclusion 
of Conjecture \ref{RIconj} holds. 
\label{RIpart}
\end{theor} 
\beginProof 
We first record the following elementary construction. Let $R$ and $T$ be two 
commutative rings and suppose that we are given a ring homomorphism $\phi:R\to T$. 
Suppose furthermore that $T$ is free as an $R$-module and let 
$t_1,\dots,t_r$ be a basis of $T$ as an $R$-module. Then there is by definition 
an isomorphism of $R$-modules 
$$
T\simeq\oplus_{k=1}^r R
$$
and the $T$-module structure of $T$ is described by a morphism or $R$-algebras 
$\psi:T\to\Mat_{r\times r}(R)$ on the right-hand side of this isomorphism. Here 
$\Mat_{r\times r}(R)$ is the ring of $r\times r$ matrices with coefficients in $R$. In particular, 
if $M$ is an $R$-module, then there is an isomorphism of $R$-modules
$$
M\otimes_R T\simeq \oplus_{k=1}^r M
$$
and the $T$-module structure of $M\otimes_R T$ is again described by $\psi$ via 
the natural action of $\Mat_{r\times r}(R)$ on $\oplus_{k=1}^r M$.

Recall that we now suppose that all the assumptions of Conjecture \ref{RIconj} 
are satisfied and that $K$ is an abelian extension of $\mQ$. 
 Let $f=f_K$ be the conductor of $K$ (in the sense of class field theory). 
We may replace wrog $D$ be a finite extension and so we may also suppose that $D$ contains some primitive $f_K$-th root of unity. Then by class field theory, there exists an embedding
$\rho:\CO_K\hookrightarrow\CO_{\mQ(\mu_f)}$ and by assumption 
there is an embedding $\lambda:\CO_{\mQ(\mu_f)}\hookrightarrow D$. We also see that every embedding of 
$K$ (resp. $\mQ(\mu_f)$) into $\mC$  factors through an embedding of $\Frac(D)$ 
into $\mC$ and similarly every embedding of 
$K$  into $\mC$  factors through an embedding of $\mQ(\mu_f)$ 
into $\mC$.

Now notice that the ring $\CO_{\mQ(\mu_f)}$ is a free module over 
$\CO_K$ via $\rho$. Indeed, $\CO_{\mQ(\mu_f)}$ is generated by a  primitive root of 
unity $z$ as an $\CO_K$-algebra. The minimal polynomial $P(X)$ of $z$ over $K$ 
divides $X^n-1$ and hence by Gauss's lemma, we have 
$P(X)\in\CO_K[X]$ and $P(X)$ is a prime element in $\CO_K[X]$. Hence there is a surjection $\CO_K[X]/(P(X))\to\CO_{\mQ(\mu_f)}$, which is also injective since $\CO_K[X]/(P(X))$ is a domain and $\CO_K[X]/(P(X))\otimes\mQ\simeq \mQ(\mn)$. Thus the elements 
$1,X,\dots X^{\deg(P)-1}$ form a basis for $\CO_{\mQ(\mu_f)}$ over $\CO_K$ via that isomorphism. 

So choose a basis $b_1,\dots b_r$ of $\CO_{\mQ(\mu_f)}$ over $\CO_K$.  
The $\CO_{\mQ(\mn)}$-module structure of $\CO_{\mQ(\mn)}$ viewed as an 
$\CO_K$-module is then described by a morphism of 
$\CO_K$-algebras $\psi:\CO_{\mQ(\mn)}\to\Mat_{r\times r}(\CO_K)$ 
(see the above elementary construction).
We let $\CB:=\times_{j=1,Y}^r\CA$ be the fibre product of $\CA$, $r$-times with itself over $Y$. The abelian scheme 
$\CB$ carries an action of $\CO_{\mQ(\mu_f)}$ via $\psi$ and we have an isomorphism 
of $\CO_{\mQ(\mu_f)}$-modules 
$$
H^1_\Db(\CB/Y)\simeq H^1_\Db(\CA/Y)\otimes_{\CO_K}\CO_{\mQ(\mu_f)}
$$
Recall that the conductor of a finite abelian extension of $\mQ$ has the same support as its discriminant. Thus by Lemma \ref{dclem}, there are decompositions into direct sums of $\CO_Y$-modules
\begin{equation}
H^1_\Db(\CB/Y)\simeq\bigoplus_{\sigma:\CO_{\mQ(\mu_f)}\hookrightarrow D}H^1_\Db(\CB/Y)_\sigma
\label{dircom}
\end{equation}
and 
$$
H^1_\Db(\CA/Y)\simeq\bigoplus_{\tau:\CO_K\hookrightarrow D}H^1_\Db(\CA/Y)_\tau.
$$
There is a natural compatibility
\begin{equation}
\bigoplus_{\sigma:\CO_{\mQ(\mu_f)}\hookrightarrow D,\,\sigma|_{\CO_K}=\tau}H^1_\Db(\CB/Y)_\sigma
\simeq \bigoplus_{j=1}^{[\mQ(\mu_f):K]}H^1_\Db(\CA/Y)_\tau
\label{findircom}
\end{equation}
where $\CO_{\mQ(\mu_f)}$ (resp. $\CO_K$) acts on $H^1_\Db(\CB/Y)_\sigma$ 
(resp. $H^1_\Db(\CA/Y)_\tau)$ via $\sigma$ (resp. $\tau$). 

Now choose a $\mu_f(\mC)$-invariant K\"ahler fibration $\kappa$ on $\CB(\mC)$ associated with a relatively ample line bundle on $\CA(\mC)$. See remark \ref{remcyc} for this. 
 
The character $\chi$ of $\Gal(K|\mQ)$ induces by composition a character $\Gal(\mQ(\mu_f)|\mQ)\to\mC$, which we shall also 
refer to as $\chi$. Choose a extension of the embedding $\iota:K\hookrightarrow D$ 
to $\mQ(\mu_f)$ and also refer to it as $\iota.$

Applying \refeq{BRR} to $\CB$ and $\kappa$ and using \refeq{findircom}, 
we obtain 
\begin{eqnarray}
&&\sum_{\sigma\in\Gal(\mQ(\mn)|\mQ)}\ari{\ch}^{[l]}(\mtr{H}^1_\Db(\CB/Y)_{\iota\circ\sigma})\chi(\sigma)
\nonumber\\&=&\nonumber-\Big[2{L'(\chi_{\rm Prim},1-l)\over L(\chi_{\rm Prim},1-l)}+\CH_{l-1}\Big]\sum_{\sigma\in\Gal(\mQ(\mn)|\mQ)}\ch^{[l-1]}(H^{1,0}(\CB/Y)_{\iota\circ\sigma})\chi(\sigma)\\
&=&
-[\mQ(\mu_f):K]\Big[2{L'(\chi_{\rm Prim},1-l)\over L(\chi_{\rm Prim},1-l)}+\CH_{l-1}\Big]\sum_{\tau\in\Gal(K|\mQ)}\ch^{[l-1]}(H^{1,0}(\CA/Y)_{\iota\circ\tau})\chi(\tau)
\label{BRR2}
\end{eqnarray}
Here $\mtr{H}^1_\Db(\CA/Y)_{\iota\circ\sigma}$ is by assumption equipped with the 
$L_2$-metric induced by the K\"ahler fibration structure $\kappa$. By Proposition 
\ref{RIconjindep} and \refeq{findircom}, the element
$$\sum_{\sigma\in\Gal(\mQ(\mn)|\mQ)}\ari{\ch}^{[l]}(\mtr{H}^1_\Db(\CB/Y)_{\iota\circ\sigma})\chi(\sigma)=
\sum_{\tau\in\Gal(K|\mQ)}\ari{\ch}^{[l]}(\mtr{H}^1_\Db(\CB/Y)_{\iota\circ\tau})\chi(\tau)
$$ does not change if we replace $\kappa$ by the K\"ahler fibration structure 
$\times_{j=1}^r\nu_f$. Hence we have
$$
\sum_{\sigma\in\Gal(\mQ(\mn)|\mQ)}\ari{\ch}^{[l]}(\mtr{H}^1_\Db(\CB/Y)_{\iota\circ\sigma})\chi(\sigma)=
[\mQ(\mu_f):K]\sum_{\tau\in\Gal(K|\mQ)}\ari{\ch}^{[l]}(\mtr{H}^1_\Db(\CA/Y)_{\iota\circ\tau})\chi(\tau)
$$
which concludes the proof. 
\endProof

\begin{rem}\rm For $l=1$, Theorem \ref{RIpart} proves a weak form of the 
conjecture of Gross-Deligne for certain linear combinations of Hodge structures 
cut out in the cohomology of $\CA(\mC)$. See \cite{MR-Periods} and also 
\cite{Bourbaki} for details. 
A different approach to this special case is described in the paper \cite{Fresan-Periods} which relies 
on a deep result of Saito and Terasoma (see \cite{ST-Det}). It would be very interesting if 
Fr\'esan's approach \cite{Fresan-Periods} could be generalised to include the case $l>1$.\end{rem}

\begin{comp}\rm The proof of Theorem \ref{RIpart} shows that under the assumptions of the Theorem \ref{RIpart}, the equality 
\begin{align*}
\sum_{\tau\in\Gal(K|\mQ)}\ari{\ch}^{[l]}(\mtr{H}^1_\Db(\CA/&Y)_{\iota\circ\tau})\chi(\tau)
\nonumber\\ =&-\Big[2\,{L'(\chi,1-l)\over L(\chi,1-l)}+\CH_{l-1}\Big]\sum_{\tau\in\Gal(K|\mQ)}\ch^{[l-1]}(H^{1,0}(\CA/Y)_{\iota\circ\tau})\chi(\tau)\nonumber
\end{align*}
actually holds in $\ari{\CHOW}^l_{\mQ(\mu_f(\chi)}(Y)$ (and not just in 
$\ari{\CHOW}^l_{\bar\mQ}(Y)$).
\end{comp}

\begin{cor} Let the assumptions of Theorem \ref{RIpart} hold. Then we have 
$$
\sum_{\tau\in\Gal(K|\mQ)}{\ch}^{[l]}({H}^1_\Db(\CA/Y)_{\iota\circ\tau})\chi(\tau)=0
$$
in $\CHOW^l(Y)_{\bar\mQ}$ for any character $\chi$ of $\Gal(K|\mQ)$ of the same parity as $l$ and any embedding 
$\iota:\CO_K\hookrightarrow D$.
\label{RIpartcor}
\end{cor}

Again, it makes sense to ask whether Corollary \ref{RIpartcor} might hold in 
a more general situation. This leads to the purely geometric

\begin{conj}
Let the assumptions of Conjecture \ref{RIconj} hold. Then $$
\sum_{\tau\in\Gal(K|\mQ)}{\ch}^{[l]}({H}^1_\Db(\CA/Y)_{\iota\circ\tau})\chi(\tau)=0
$$
in $\CHOW^l(Y)_{\bar\mQ}$ for any character $\chi$ of $\Gal(K|\mQ)$ of the same parity as $l$ and any embedding 
$\iota:\CO_K\hookrightarrow D$.
\end{conj}

See \cite[Prop. 3]{MR-Order} for more about this conjecture in a slightly more restrictive setting.

We now indulge in some wilder speculation. The fact that the formula in Conjecture 
\ref{RIconj} looks very 'motivic' suggests the following vague conjecture, which 
seems to be a good computational thumb rule in many examples. 

\begin{vconj}
Suppose that $K$ is a finite Galois extension of $\mQ$. Suppose that 
all the embedding of $K$ into $\mC$ factor through an embedding 
of $\Frac(D)$ into $\mC$. Suppose also that the discriminant of $K$ is invertible in $D$ and that 
$D$ is a localisation of the ring of integers of a number field. 

Let $f:\CM\to Y$ be a 'log smooth relative motive' over Y and suppose that 
we are given an embedding of rings $\CO_K\hookrightarrow\End_Y(\CM).$ 

Let 
$\chi:\Gal(K|\mQ)\to\mC$ be an irreducible Artin character. Let $l\geq 1$. 
Suppose that $\chi$ and $l$ have the same parity (hence $L(\chi,1-l)\not=0$ by Remark \ref{rem_parity}). 
 
Then there exists a 'polarisation' on $\CM/Y$, which is compatible with 
the action of $\CO_K$ in some sense and such that for any embedding $\iota:\CO_K\hookrightarrow D$ we have:
\begin{align*}
\sum_{\tau\in\Gal(K|\mQ)}\ari{\ch}^{[l]}&(\mtr{H}^k_\Db(X/Y)_{\iota\circ\tau}(\log))\chi(\tau)
\nonumber\\ =&-\Big[2\,{L'(\chi,1-l)\over L(\chi,1-l)}+\CH_{l-1}\Big]\sum_{p+q=k}\sum_{\tau\in\Gal(K|\mQ)}p\cdot\ch^{[l-1]}(H^{p,q}(X/Y)_{\iota\circ\tau}(\log))\chi(\tau)\nonumber
\end{align*}
in $\ari{\CHOW}^l_{K(\chi)}(Y)(\log).$ 
\label{vRIconj}
 \end{vconj}

Here ${H}^k_\Db(X/Y)(\log)$ and $H^{p,q}(X/Y)(\log)$ refers to logarithmic cohomology and 
the metric on ${H}^k_\Db(X/Y)(\log)$ is induced by the polarisation, which is general mildly singular. The ring $\ari{\CHOW}^l_{K(\chi)}(Y)(\log)$ is a generalised arithmetic 
intersection ring, as in \cite{BKK}.
Note that in this vague conjecture, if $\CM$ is smooth over $Y$, then 
one may remove the '$(\log)$' symbols from the formula.

In particular, this 'conjecture' should apply to generically abelian semiabelian schemes, where it 
should be possible to make a precise conjecture, extending Conjecture \ref{RIconj}. We refrain from trying to do this 
here because the generalised arithmetic 
intersection theory that would be necessary for this has apparently not yet been fully defined. 
In some of the examples drawn from the literature that we shall consider in section \ref{secex} below, the relevant articles produce generalised arithmetic 
intersection theories tailor-made for the 
geometric situation under consideration.

\section{Examples}
\label{secex}

%
We shall now show that various formulae proven in the literature are formally compatible 
with Conjecture \ref{RIconj} and in some cases are partially consequences of Theorem \ref{RIpart}. We use the notations of Conjecture \ref{RIconj}. We 
shall make the assumption that in each example considered below the polarisation has been chosen 
in such a way that for each $\tau\in\Gal(K|\mQ)$ we have an isometric isomorphism 
$$
f_*(\Omega_{\CA/Y})_{\iota\circ\tau}\simeq \R^1 f_*(\CO_\CA)_{\iota\circ c \circ\tau}^\vee
$$ 
This is a compatibility with duality that is often verified in practice. 

{\bf The formula of Colmez.} (see \cite{Colmez-Per} and \cite{Colmez-Hauteur}) Suppose that $Y=\Spec\,D$ and that $K$ is a 
CM field of degree $2\cdot\dim(\CA/Y).$ Suppose that $K$ is an abelian extension of $\mQ$. Let $\Phi:\Hom(K,D)\to\{0,1\}$ be the 
associated CM type. By definition, 
$$
\Phi(\iota\circ\tau)=\rk(H^{1,0}(X/Y)_{\iota\circ c\circ\tau}).
$$
We identify $\Phi$ with a function $\Gal(K|\mQ)\to\{0,1\}$ via $\iota$. 
From now until the end of the computation, we shall drop the embedding $\iota$ from the notation. Theorem \ref{RIpart} gives the equality:
\begin{align*}
\sum_{\tau\in\Gal(K|\mQ)}\ac1(\mtr{H}^1_\Db(\CA/Y)_\tau)\chi(\tau)&=-2\,{L'(\chi,0)\over L(\chi,0)}\sum_{\tau\in\Gal(K|\mQ)}\Phi(\tau)\chi(c\circ \tau)
\\
&=
2\,{L'(\chi,0)\over L(\chi,0)}\sum_{\tau\in\Gal(K|\mQ)}\Phi(\tau)\chi(\tau)
\end{align*}
in $\ari{\rm CH}^1_{\bar\mQ}(D)$ 
for any odd one-dimensional character. By assumption, we have
$$
\ac1(\mtr{H}^1_\Db(\CA/Y)_\tau)=-\ac1(\mtr{H}^1_\Db(\CA/Y)_{c\circ\tau}).
$$
so that  
for even characters $\chi$, we have
$$
\sum_{\tau\in\Gal(K|\mQ)}\ac1(\mtr{H}^1_\Db(\CA/Y)_\tau)\chi(\tau)=0.
$$
We recall the definition of the scalar product
$$
\langle f,g\rangle:={1\over [K:\mQ]}\sum_{\tau\in\Gal(K|\mQ)}f(\tau)\overline{g(\tau)}
$$
and of the convolution product
$$
(f\ast g)(\sigma):={1\over [K:\mQ]}\sum_{\tau\in\Gal(K|\mQ)}g(\tau)f(\tau^{-1}\sigma)
$$
of two functions $f,g:\Gal(K|\mQ)\to\mC$. 
Recall that if $h$ is a one-dimensional character then we have 
$$
\langle f\ast g,h\rangle=\langle f,h\rangle\cdot\langle g,h\rangle
$$
for any two functions $f,g:\Gal(K|\mQ)\to\mC.$

Define the function $\Phi^\vee:\Gal(K|\mQ)\to\{0,1\}$ by the formula 
$\Phi^\vee(\tau):=\Phi(\tau^{-1})$. 
Using the fact that the one-dimensional characters $\Gal(K|\mQ)\to\mC$ form an orthogonal basis
 of  the space of complex valued functions on $\Gal(K|\mQ)$ we get that
\begin{equation}
\sum_{\tau}\ac1(\mtr{H}^1_\Db(\CA/Y)_\tau)={1\over [K:\mQ]}
\sum_{\chi\,{\rm odd}}\bar\chi(\tau)\Big[2\,{L'(\chi,0)\over L(\chi,0)}\sum_{\sigma\in\Gal(K|\mQ)}\Phi(\sigma)\chi(\sigma)\Big]
\label{IDKneq}
\end{equation}
 From \refeq{IDKneq} we obtain the equality 
\begin{eqnarray}
&&\ac1(f_*(\Omega_{\CA/Y}))=\sum_{\tau\in\Gal(K|\mQ)}\ac1(\mtr{H}^1_\Db(\CA/Y)_\tau)\Phi(c\circ \tau)\nonumber=
-\sum_{\tau\in\Gal(K|\mQ)}\ac1(\mtr{H}^1_\Db(\CA/Y)_\tau)\Phi(\tau)\nonumber\\&=&
-{1\over [K:\mQ]}\sum_{\chi\,\odd}2\,{L'(\chi,0)\over L(\chi,0)}\Big(\sum_{\tau\in\Gal(K|\mQ)}\Phi(\tau)\bar\chi(\tau)\Big)\Big(\sum_{\sigma\in\Gal(K|\mQ)}\Phi(\sigma)\chi(\sigma)\Big)\nonumber\\
&=&
-{1\over [K:\mQ]}\sum_{\chi\,\odd}2\,{L'(\chi,0)\over L(\chi,0)}\Big(\sum_{\tau\in\Gal(K|\mQ)}\Phi(\tau)\bar\chi(\tau)\Big)\Big(\sum_{\sigma\in\Gal(K|\mQ)}\Phi^\vee(\sigma)\bar\chi(\sigma)\Big)\nonumber\\
&=&
-{[K:\mQ]}\cdot \sum_{\chi\,\odd}2\,{L'(\chi,0)\over L(\chi,0)}\langle \Phi\ast\Phi^\vee,\chi\rangle.
\label{colgen}
\end{eqnarray}
The formula \refeq{colgen} implies the formula of Colmez (see \cite[Conjecture 3]{Colmez-Hauteur} and the discussion after the statement) up to a term of the form 
$$
\sum_{p|D_K}r_p\log(p)
$$
where $r_p\in\bar\mQ$.

{\bf The formula of Bost and Kühn.} (see \cite{Kuehn} and also an unpublished manuscript by J.-B. Bost) In that case, $\CA/Y$ is an elliptic scheme and 
$K=\mQ$. Applying Theorem \ref{RIpart}, we obtain

\begin{eqnarray*}
\ari{\ch}^{[2]}(\bar H^1(\CA/Y))=-\Big[2\,{\zeta_\mQ'(-1)\over\zeta_\mQ(-1)}+1\Big]\cdot{\rm c}^1(H^{1,0}(\CA/Y))
\end{eqnarray*}
ie
\begin{eqnarray}
\ari{\rm c}^1(\bar{\omega})^2=-\Big[2\,{\zeta_\mQ'(-1)\over\zeta_\mQ(-1)}+1\Big]\cdot{\rm c}^1(\omega)
\label{BKfor}
\end{eqnarray}
where $\bar\omega$ is the Hodge bundle of $\CA$ ({\it i.e.} the restriction of the sheaf of 
differentials of $\CA/Y$ by the unit section) endowed with the Peterson metric. 
Note that equality \refeq{BKfor} is of little interest because if $\CA_\mC$ has non zero Kodaira-Spencer class then $Y$ cannot be chosen to be proper over $D$ (this follows from the structure of the moduli spaces of elliptic curves) so that one always has 
${\rm c}^1(\omega)=0$. The formula of Bost and Kühn has the same shape as \refeq{BKfor} but is valid for some generically abelian semiabelian schemes over $Y$ 
(for which ${\rm c}^1(\omega)\not=0$). It 
allows mild singularities and can thus be understood as a 'special case' of the 
vague conjecture \ref{vRIconj}.

{\bf Families of abelian surfaces with complex multiplication by a quadratic imaginary extension of 
$\mQ$. The formula of Kudla, Rapoport and Yang.} (see \cite[T. 1.05]{KRY-Mod})

In that case, $\dim(\CA/Y)=2$ and $K$ is a quadratic imaginary extension of $\mQ$. 
In particular the group $\Gal(K|\mQ)$ has precisely one non-trivial character $\chi$ and this character is odd. 
Theorem \ref{RIpart} gives 
\begin{eqnarray}
&&\sum_{\tau\in\Gal(K|\mQ)}\ari{\ch}^{[1]}(\mtr{H}^1_\Db(\CA/Y)_\tau)\chi(\tau)=-2\,{L'(\chi,0)\over L(\chi,0)}\sum_{\tau\in\Gal(K|\mQ)}\rk(H^{1,0}(\CA/Y)_\tau)\chi(\tau)
\label{KRYeq}
\end{eqnarray}
Write $\bar\tau:=c\circ \tau$ in the following computations. From now on until the end of the computation, we shall drop the embedding $\iota$ from the notation. 
We have a decomposition
$$
f_*(\Omega_{\CA/Y})\simeq f_*(\Omega_{\CA/Y})_\tau\oplus f_*(\Omega_{\CA/Y})_{\bar\tau}
$$
and 
$$
\R^1 f_*(\CO_\CA)\simeq  
(f_*(\Omega_{\CA/Y})_\tau)^\vee\oplus (f_*(\Omega_{\CA/Y})_{\bar\tau})^\vee
$$
so that we may rewrite \refeq{KRYeq} as 
\begin{eqnarray*}
&&2\Big(\ac1(f_*(\bar\Omega_{\CA/Y})_\tau)-\ac1(f_*(\bar\Omega_{\CA/Y})_{\bar\tau})\Big)
=-2\,{L'(\chi,0)\over L(\chi,0)}\sum_{\tau\in\Gal(K|\mQ)}\rk(f_*(\Omega_{\CA/Y})_\tau)\chi(\tau)
\end{eqnarray*}
Squaring the preceding equality, we see that 
\begin{eqnarray*}
&&\Big(\ac1(f_*(\bar\Omega_{\CA/Y})_\tau)-\ac1(f_*(\bar\Omega_{\CA/Y})_{\bar\tau})\Big)^2\\&&=
\ac1(f_*(\bar\Omega_{\CA/Y})_\tau)^2+\ac1(f_*(\bar\Omega_{\CA/Y})_{\bar\tau})^2-
2\cdot \ac1(f_*(\bar\Omega_{\CA/Y})_\tau)\cdot \ac1(f_*(\bar\Omega_{\CA/Y})_{\bar\tau})=0
\end{eqnarray*}
so that 
\begin{eqnarray*}
\ac1(f_*(\bar\Omega_{\CA/Y})_\tau)^2+\ac1(f_*(\bar\Omega_{\CA/Y})_{\bar\tau})^2=
2\cdot \ac1(f_*(\bar\Omega_{\CA/Y})_\tau)\cdot \ac1(f_*(\bar\Omega_{\CA/Y})_{\bar\tau}).
\end{eqnarray*}
Now since $\CA/Y$ can also be viewed as carrying an action of $\mQ=\mQ(\mu_2)$, Theorem \ref{RIpart} also gives
\begin{eqnarray*}
\ari{\ch}^{[2]}(\bar H^1(\CA/Y))=-\Big[2\,{\zeta_\mQ'(-1)\over\zeta_\mQ(-1)}+1\Big]\cdot{\rm c}^1(H^{1,0}(\CA/Y))=-\Big[2\,{\zeta_\mQ'(-1)\over\zeta_\mQ(-1)}+1\Big]\cdot{\rm c}^1(f_*(\Omega_{\CA/Y}))
\end{eqnarray*}
where now
$$
\ari{\ch}^{[2]}(\bar H^1(\CA/Y))=\ac1(f_*(\bar\Omega_{\CA/Y})_\tau)^2+\ac1(f_*(\bar\Omega_{\CA/Y})_{\bar\tau})^2
$$
so that 
$$
\ac1(f_*(\bar\Omega_{\CA/Y})^2=
(\ac1(f_*(\bar\Omega_{\CA/Y})_\tau)+\ac1(f_*(\bar\Omega_{\CA/Y})_{\bar\tau}))^2=
2\cdot \ac1(f_*(\bar\Omega_{\CA/Y})_\tau)^2+2\cdot \ac1(f_*(\bar\Omega_{\CA/Y})_{\bar\tau})^2
$$
and
\begin{equation}
\ac1(f_*(\bar\Omega_{\CA/Y})^2=-2\cdot\Big[2\,{\zeta_\mQ'(-1)\over\zeta_\mQ(-1)}+1\Big]\cdot{\rm c}^1(f_*(\Omega_{\CA/Y})).
\label{finKRY}
\end{equation}
Formula \refeq{finKRY} implies the formula \cite[Th. 1.0.5]{KRY-Mod} up to a factor of the form $$
\sum_{p|D_K}r_p\log(p)
$$
where $r_p\in\bar\mQ$.

{\bf Families of abelian surfaces with an action by a real quadratic extension of $\mQ$. The formula of 
Bruiner, Burgos and K\"uhn.} (see \cite[Th. B]{BBK-Bor}) In this case $\dim(\CA/Y)=2$ and $K$ is a real quadratic extension of $\mQ$. All the one-dimensional characters of $$\Gal(K|\mQ)=\{\Id,\tau_0\}$$ are even and there is only one non-trivial one-dimensional character $\chi_0$. We again drop the embedding 
$\iota:K\hookrightarrow D$ from the notations. 
For any one-dimensional character, Theorem \ref{RIpart} gives:
\begin{eqnarray*}
&&\sum_{\tau\in\Gal(K|\mQ)}\ari{\ch}^{[2]}(\mtr{H}^1_\Db(\CA/Y)_\tau)\chi(\tau)=-\Big[2\,{L'(\chi,-1)\over L(\chi,-1)}+\CH_{l-1}\Big]\sum_{\tau\in\Gal(K|\mQ)}{\rm c}^1(H^{1,0}(\CA/Y)_\tau)\chi(\tau)
\end{eqnarray*}
By assumption, this translates to
\begin{eqnarray*}
&&2\sum_{\tau\in\Gal(K|\mQ)}\ari{\ch}^{[2]}((f_*(\mtr{\Omega}_{\CA/Y}))_\tau)\chi(\tau)=-\Big[2\,{L'(\chi,-1)\over L(\chi,-1)}+1\Big]\sum_{\tau\in\Gal(K|\mQ)}{\rm c}^1((f_*(\Omega_{\CA/Y}))_\tau)\chi(\tau)
\end{eqnarray*}
Specialising this to each character, we obtain:
\begin{eqnarray*}
&&\ac1((f_*(\Omega_{\CA/Y}))_{\Id})^2+\ac1((f_*(\Omega_{\CA/Y}))_{\tau_0})^2=-\Big[2\,{\zeta_\mQ'(-1)\over \zeta_\mQ(-1)}+1\Big]{\rm c}^1(f_*(\Omega_{\CA/Y}))
\end{eqnarray*}
and
\begin{multline*}
\ac1((f_*(\Omega_{\CA/Y}))_{\Id})^2-\ac1((f_*(\Omega_{\CA/Y}))_{\tau_0})^2
\\
=-\Big[2\,{L'(\chi_0,-1)\over L(\chi_0,-1)}+1\Big]\Big({\rm c}^1((f_*(\Omega_{\CA/Y}))_{\Id})-{\rm c}^1((f_*(\Omega_{\CA/Y}))_{\tau_0})\Big).
\end{multline*}
Note that this implies that
$${\rm c}^1((f_*(\Omega_{\CA/Y}))_{\tau_0})^2={\rm c}^1((f_*(\Omega_{\CA/Y}))_{\Id})^2=0.$$
Now we may compute
\begin{eqnarray*}
&&
\ac1(f_*(\mtr{\Omega}_{\CA/Y}))^3=\Big(\ac1((f_*(\Omega_{\CA/Y}))_{\Id})+\ac1((f_*(\Omega_{\CA/Y}))_{\tau_0})\Big)^3\\&=&
\ac1((f_*(\Omega_{\CA/Y}))_{\Id})^3+\ac1((f_*(\Omega_{\CA/Y}))_{\tau_0})^3+
3\cdot \ac1((f_*(\Omega_{\CA/Y}))_{\tau_0})\cdot\ac1((f_*(\Omega_{\CA/Y}))_{\Id})^2\\&+&3\cdot \ac1((f_*(\Omega_{\CA/Y}))_{\Id})\cdot\ac1((f_*(\Omega_{\CA/Y}))_{\tau_0})^2\\
&=&
-\Big({\rm c}^1((f_*(\Omega_{\CA/Y}))_{\Id})+3\cdot{\rm c}^1((f_*(\Omega_{\CA/Y}))_{\tau_0})\Big)\cdot
{1\over 2}\Big[[2\,{\zeta_\mQ'(-1)\over \zeta_\mQ(-1)}+1]{\rm c}^1(f_*(\Omega_{\CA/Y}))\\&+&
[2\,{L'(\chi_0,-1)\over L(\chi_0,-1)}+1]({\rm c}^1((f_*(\Omega_{\CA/Y}))_{\Id})-{\rm c}^1((f_*(\Omega_{\CA/Y}))_{\tau_0}))\Big]\\
&-&
\Big({\rm c}^1((f_*(\Omega_{\CA/Y}))_{\tau_0})+3\cdot{\rm c}^1((f_*(\Omega_{\CA/Y}))_{\Id})\Big)\cdot
{1\over 2}\Big[[2\,{\zeta_\mQ'(-1)\over \zeta_\mQ(-1)}+1]{\rm c}^1(f_*(\Omega_{\CA/Y}))\\&-&
[2\,{L'(\chi_0,-1)\over L(\chi_0,-1)}+1]({\rm c}^1((f_*(\Omega_{\CA/Y}))_{\Id})-{\rm c}^1((f_*(\Omega_{\CA/Y}))_{\tau_0}))\Big]\\
&=&
-\Big(2\cdot[2\,{\zeta_\mQ'(-1)\over \zeta_\mQ(-1)}+1]+[2\,{L'(\chi_0,-1)\over L(\chi_0,-1)}+1]\Big)\cdot{\rm c}^1(\Omega_{\CA/Y})^2
\end{eqnarray*}
This may be rewritten in terms of the zeta function of $K$. Recall that we have
$$
\zeta_K(s)=\zeta_\mQ(s)L(\chi_0,s).
$$
We finally obtain the equality
\begin{equation}
\ac1(f_*(\mtr{\Omega}_{\CA/Y}))^3=-\Big(2\,{\zeta_\mQ'(-1)\over \zeta_\mQ(-1)}+2\,{\zeta_K'(-1)\over \zeta_K(-1)}+3\Big)\cdot{\rm c}^1(\Omega_{\CA/Y})^2.
\label{BBKfor}
\end{equation}
which should be compared with \cite[Th. B]{BBK-Bor}. 
As for the formula of Bost and K\"uhn, equality \refeq{BBKfor} is not very interesting because $\CA$ is not allowed to be semiabelian. The formula in  \cite[Th. B]{BBK-Bor} has the same shape 
as \refeq{BBKfor} but allows semiabelian schemes and allows the metric to have mild singularities. It can thus again be understood as a 'special case' of the vague conjecture 
\ref{vRIconj}.

\end{document}